\documentclass[10pt,twocolumn]{article}

\usepackage[utf8]{inputenc}
\usepackage[T1]{fontenc}
\usepackage{lmodern}
\usepackage[margin=1in]{geometry}
\usepackage{setspace}
\usepackage{authblk}
\usepackage[hidelinks]{hyperref}
\usepackage{titling}
\usepackage{abstract}
\usepackage[title]{appendix}
\usepackage{placeins}
\usepackage{cuted}
\usepackage[backend=biber,style=numeric,maxnames=5, minnames=3]{biblatex}

\title{Hybrid Quantum-Classical Optimization for \\ Multi-Objective Supply Chain Logistics}

\author{%
  Raoul Heese\textsuperscript{1,*},
  Timoth\'ee Leleu\textsuperscript{2,3},
  Sam Reifenstein\textsuperscript{2}, \\
  Christian Nietner\textsuperscript{1}, and
  Yoshihisa Yamamoto\textsuperscript{2}
}

\date{%
  {\small
  \textsuperscript{1}NTT Data Germany\\
  \textsuperscript{2}NTT Research\\
  \textsuperscript{3}Stanford University\\[0.5em]
  \textsuperscript{*}corresponding author: raoul.heese@nttdata.com}
  }

\usepackage{float,subcaption}
\usepackage{csquotes}
\usepackage{mathtools,amsmath,amsfonts,amssymb,optidef,braket}
\usepackage{cuted}
\usepackage{algorithm,algpseudocode}
\usepackage{mdframed}
\usepackage{siunitx}
\usepackage{textcomp}
\usepackage{enumitem}
\usepackage{booktabs}
\usepackage{pifont}
\usepackage[dvipsnames]{xcolor}
\usepackage{tikz}
\usepackage{pgfplots}
\usetikzlibrary{calc,decorations.pathreplacing,patterns,shapes.geometric,arrows.meta}
\usepgfplotslibrary{groupplots}
\pgfplotsset{compat=newest}
\usepackage{pgfplotstable}
\usepackage[mode=buildnew]{standalone}
\usepackage{newpxtext}
\usepackage{ifthen}
\usepackage{url}
\usepackage[acronym]{glossaries}
\usepackage[hidelinks]{hyperref}
\usepackage[nameinlink,capitalise]{cleveref} 

\makeatletter
\newcommand{\ie}{i.e.\@ifnextchar.{\!\@gobble}{}}
\newcommand{\eg}{e.g.\@ifnextchar.{\!\@gobble}{}}
\newcommand{\etc}{etc\@ifnextchar.{}{.\@}}
\makeatother

\DeclareMathOperator*{\argmin}{arg\,min}
\DeclareMathOperator{\sign}{sign}

\newacronym{OR}{OR}{Operations Research}
\newacronym{PBS}{PBS}{Product Breakdown Structure}
\newacronym{KPI}{KPI}{Key Performance Indicator}
\newacronym{QUBO}{QUBO}{Quadratic Unconstrained Binary Optimization Problem}
\newacronym{QCBO}{QCBO}{Quadratic Constrained Binary Optimization Problem}
\newacronym{VQE}{VQE}{Variational Quantum Eigensolver}
\newacronym{QAOA}{QAOA}{Quantum Approximate Optimization Algorithm}
\newacronym{ISG}{ISG}{Informed Solution Generator}
\newacronym{ISF}{ISF}{Informed Solution Fixer}
\newacronym{ISI}{ISI}{Informed Solution Improver}
\newacronym{IQTS}{IQTS}{Informed Quantum-Enhanced Tree Solver}
\newacronym{HBS}{HBS}{Hybrid Quantum-Classical Bilevel Solver}
\newacronym{DAS}{DAS}{Dynamic Anisotropic Smoothing}
\newacronym{BP}{BP}{Belief propagation}
\newacronym{CACm}{CACm}{Chaotic Amplitude Control with Momentum}
\newacronym{IBP}{IBP}{Iterative Belief Propagation}
\newacronym{ODE}{ODE}{Ordinary Differential Equation}
\newacronym{CIM}{CIM}{Coherent Ising Machine}
\newacronym{TFLN}{TFLN}{Thin Film Lithium Niobate}
\newacronym{NISQ}{NISQ}{Noisy Intermediate-Scale Quantum}
\newacronym{SA}{SA}{Simulated Annealing}
\newacronym{PT}{PT}{Parallel Tempering}
\newacronym{VQA}{VQA}{Variational Quantum Algorithm}

\newlist{enumeratedense}{enumerate}{3}
\setlist[enumeratedense]{noitemsep,nolistsep}
\setlist[enumeratedense,1]{label=\arabic*.}
\setlist[enumeratedense,2]{label=\alph*.}
\newlist{itemizedense}{itemize}{3}
\setlist[itemizedense]{noitemsep,nolistsep}
\setlist[itemizedense,1]{label=\textbullet}
\newlist{experimentlist}{enumerate}{1}
\setlist[experimentlist]{noitemsep,nolistsep}
\setlist[experimentlist,1]{label=E\arabic*.,ref=E\arabic*,leftmargin=1.2cm}
\crefname{experimentlisti}{experiment}{experiments}

\newboolean{ejorflag}
\newboolean{arxivflag}
\newenvironment{condstrip}{\ifthenelse{\boolean{arxivflag}}{\begin{strip}}{}}
{\ifthenelse{\boolean{arxivflag}}{\end{strip}}{}}
\newcommand{\condclearpage}{\ifthenelse{\boolean{arxivflag}}{\clearpage}{}}
\setboolean{ejorflag}{false}
\setboolean{arxivflag}{true}

\begin{document}      

\maketitle

\begin{strip}
\\[-2cm]
\begin{abstract}
A multi-objective logistics optimization problem from a real-world supply chain is formulated as a Quadratic Unconstrained Binary Optimization Problem (QUBO) that minimizes cost, emissions, and delivery time, while maintaining target distributions of supplier workshare. The model incorporates realistic constraints, including part dependencies, double sourcing, and multimodal transport. Two hybrid quantum-classical solvers are proposed: a structure-aware informed tree search (IQTS) and a modular bilevel framework (HBS), combining quantum subroutines with classical heuristics. Experimental results on IonQ's Aria-1 hardware demonstrate a methodology to map real-world logistics problems onto emerging combinatorial optimization-specialized hardware, yielding high-quality, Pareto-optimal solutions.
\end{abstract}
\end{strip}

\section{Introduction}

The interdependent structure of global supply chains for complex manufacturing is characterized by multi-factorial requirements that address economic, environmental, social, and political constraints. These factors require careful balancing of potentially conflicting objectives, such as cost, sustainability, and resilience~\cite{barbosapovoa2018supplychain}. The underlying combinatorial problems in global supply chain optimization are typically NP-hard, making them notoriously difficult to solve at scale. Recent research has turned to novel computing paradigms, including specialized hardware like Ising machines~\cite{mohseni2022ising} and quantum annealers~\cite{abbas2024}, to address these challenges more efficiently.

However, solving real-world optimization problems with these emerging platforms remains a considerable open problem due to current technological limitations~\cite{abbas2024}. As a result, identifying suitable applications for non-classical computing approaches is an active area of research~\cite{koch2025optimization}. In this context, it is important to explore use cases that are both sufficiently hard to solve and capable of delivering real-world impact if solved more effectively. These criteria apply to supply chain optimization, making it a promising use case candidate for further research. In this paper, a hybrid quantum-classical approach is developed to solve a real-world, multi-objective logistics optimization problem.

Our approach combines specialized modeling with hybrid quantum-classical and quantum-inspired solvers. We integrate high-performance classical optimization techniques with advanced quantum methods, allowing us to estimate the Pareto frontier of a multi-objective optimization problem. For our numerical evaluations, we use \emph{Aria-1}~\cite{chen2024benchmarking}, a state-of-the-art and commercially deployed quantum hardware trapped-ion quantum computer from IonQ. Since we focus on scalability and modularity, our findings are also applicable to a more general scope, opening up opportunities to leverage special-purpose Ising machines for complex real-world optimization problems. 

Our main contributions are as follows. We developed a multi-objective \gls{QUBO} model, aiming to optimize four \glspl{KPI} of a logistics network: carbon dioxide emissions, costs, time, and supplier target workshare fulfillment. The model integrates real-world constraints such as mandatory workshare limits, multiple transportation modes, and dependencies between parts. The formulation as a \gls{QUBO} was driven by the enabling of quantum algorithms. Next, we designed specialized solver tools---\gls{ISG}, \gls{ISF}, and \gls{ISI}---that exploit problem-specific knowledge to construct, repair, and refine solutions iteratively, improving solution quality and feasibility. We propose two hybrid quantum-classical solvers: \gls{IQTS} and \gls{HBS}. \Gls{IQTS} aims at exploiting our knowledge about the problem structure to achieve a very fast convergence. It is tightly tailored to the problem by integrating the \gls{QAOA} into an explorative tree decomposition. \Gls{HBS} is a general high-performance strategy with a highly modular architecture. It is aimed at providing computational speed, and scalability by combining quantum methods such as \gls{QAOA}, quantum-inspired methods such as \gls{CACm}, and classical methods such as \gls{IBP} and \gls{DAS} into a bilevel optimization framework. We test our algorithms on an industrially relevant logistic network provided by Airbus and are able to find Pareto-optimal solution candidates.

The remaining paper is structured a follows. In \cref{sec:related work}, we outline related work. Subsequently, in \cref{sec:use case}, we describe the use case. Our solution approaches are then presented in \cref{sec:methods}. In \cref{sec:experiments}, we show the results of our numerical evaluations and experiments. We end with a conclusion in \cref{sec:conclusions}.

\section{Related work}\label{sec:related work}

Quantum optimization in an emerging field~\cite{abbas2024} with many different application domains, mostly focusing on \glspl{QUBO}~\cite{desantis2024}. Multi-objective \glspl{QUBO} have already been considered~\cite{ayodele2022}, also with respect to different scalarisation~\cite{ayodele2023} and standardization techniques~\cite{lee2025}.

Significant effort has been dedicated in recent years to advancing hybrid quantum-classical methodologies for logistics optimization, particularly to address the computational limitations of current quantum hardware. A hybrid workflow has been proposed to tackle real-scale multi-truck vehicle routing problems through iterative quantum computations, making it feasible for \gls{NISQ} devices by segmenting complex tasks into manageable quadratic optimization instances suitable for quantum annealers~\cite{weinberg2023supply}. Similarly, quantum annealing has been used effectively for the Transport Network Design Problem (TNDP), formulated as \gls{QUBO}, providing computational benefits over classical Tabu search methods, with implications extending across various network optimization scenarios~\cite{dixit2023quantum}. Furthermore, a quantum-classical strategy called \emph{Q4RPD} has been developed, specifically addressing realistic constraints in vehicle routing, such as heterogeneous vehicle fleets, utilizing the D-Wave Leap solver~\cite{osaba2024solving}. In addition, \gls{QAOA} has been applied to the Heterogeneous Vehicle Routing Problem (HVRP), demonstrating potential for combinatorial optimization, despite facing scalability challenges related to the required number of qubits~\cite{fitzek2024applying}.

Variational approaches, such as \gls{QAOA} and advanced forms of \gls{VQE}, have been employed to solve job shop scheduling and vehicle routing problems, demonstrating the feasibility of these techniques on current quantum hardware and achieving near-optimal results within small-scale instances~\cite{amaro2022case,kurowski2023application}. In parallel, hybrid annealing frameworks have proven effective in addressing multi-agent pathfinding and complex delivery routing under realistic constraints, with iterative quantum subproblem solving enabling scalability to industry-relevant scenarios~\cite{gerlach2025hybrid}. Moreover, hybrid learning architectures integrating quantum circuits into reinforcement learning agents or neural networks have achieved competitive performance in logistics applications like fleet dispatch and backorder prediction, revealing promising intersections between quantum machine learning and supply chain analytics~\cite{correll2023quantum,jahin2023qamplifynet}.

\section{Use case} \label{sec:use case}

The 2024 \emph{Quantum Computing Challenge} from Airbus and BMW~\cite{QuantumChallenge2024} contained a task in which the goal was to optimize a logistics network for a sustainable aircraft manufacturing process within a global supply chain.

One aspect of the supply chain of the Airbus aircraft manufacturing process is to decide, by which supplier at which global production sites, the aircraft parts are produced and how they are transported between these sites for assembling. The goal is to minimize costs, production time, carbon dioxide emissions, and supplier target workshare fulfillments while ensuring certain constraints that require that the assignments are sufficiently distributed both between all suppliers and all production sites. This logistics problem is a particularly challenging aspect of the supply chain optimization and was proposed by Airbus. In the following, the Airbus supply chain problem is used to illustrate the proposed method, which is general and applicable to a broader class of optimization problems.


\section{Model} \label{sec:model}

To begin with, we first present a detailed problem description, followed by a formal definition of the properties of a problem instance. Subsequently, we provide a multi-objective \gls{QCBO} formulation and show how it can be reduced to a scalarized \gls{QUBO} form.

\subsection{Problem description}

\begin{figure*}[!tb]
    \centering
    \includegraphics[scale=1]{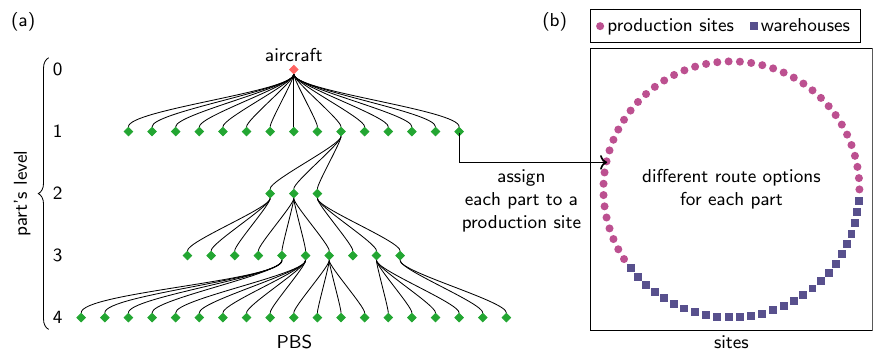}
    \caption{Problem sketch: (a) The \gls{PBS} showing dependencies of all \num{48} parts with the assembled product, an aircraft, as the root. (b) Supply chain sites comprised of \num{43} production sites and \num{28} warehouses. The shown positions do not reflect geographical locations.}
    \label{fig:problemSketch}
\end{figure*}

The supply chain logistics model we consider is based on the production of a single product. This product is assembled from a number of individual parts—48 in the case of the Airbus aircraft. The so-called \gls{PBS} shown in \cref{fig:problemSketch}a describes the dependency of the parts as a tree structure. Each node is a part and its child nodes indicate which other parts are required for its production. The parts represented by leaf nodes can therefore be produced without any dependencies. The completed aircraft is shown as the root. Each part can be assigned a level (from 0 to 4) that indicates how many production steps it takes for it to be integrated into the final product.

In total, there are \num{29} suppliers, each of which can produce a specific selection of parts at specific production sites. To produce a part at a specific site, all parts it depends on must be available at the same site. To that end, parts can be transported from one site to another. Transportation can be done by different means depending on the start and end site and the part to be transported. Options include trucks, ships, and cargo planes. Each mode of transportation is associated with carbon dioxide emissions, costs, and transportation time. A part can either be transported directly or indirectly via certain hubs, called warehouses, where it can be transshipped. Consequently, a large variety of transportation options can be selected for each part. We consider \num{43} production sites and \num{28} warehouses as sketched in \cref{fig:problemSketch}b, which are distributed globally across \num{18} different regions.

To improve the supply chain resilience, it is necessary to install double sourcing, if possible. This means that the production of each part is proportionally split between two production sites (and suppliers), the \emph{primary} site (supplier) and the \emph{secondary} site (supplier). While the two sites (primary production site and secondary production site) must be different, the respective suppliers (primary supplier and secondary supplier) can be the same. Additionally, if possible, the two sites must be located in different regions. Regardless of whether it is produced as a primary or secondary production, each part is transported to both the primary and secondary processing sites, resulting in up to four necessary transportation routes per part. The only exception are two immobile parts, which can only be produced and processed at the same site, which they will never leave.

Each part is assigned a relative cost that measures its value in relation to the other parts. This relative cost is used to ensure a shared distribution of the production process among different production sites and suppliers. For this purpose, each manufacturing site and each supplier has a predefined lower and upper limit for the cumulative costs of the products produced, known as the \emph{workshare}. In addition, each supplier has a fixed target for the workshare that should ideally be achieved.

In summary, the task is a multi-objective assignment problem in which one needs to (i) assign each part to a primary and secondary production site and supplier, and (ii) assign transportation routes to move parts between production sites (possibly via warehouses). The assignments must meet the following criteria: (i) compliance with the \gls{PBS} dependencies, (ii) adherence to workshare fulfillment constraints for production sites and suppliers, and (iii) the primary and secondary sites for each part must be distinct and, if possible, located in different regions. The goal is to jointly optimize four \glspl{KPI}: cumulative carbon dioxide emissions, costs, production time, and supplier target workshare fulfillments. The first two \glspl{KPI} are determined by the chosen transportation routes alone, while production time is affected by both the chosen transportation routes and the selected production sites and suppliers. Lastly, target workshare fulfillment is determined by the supplier choices.

\subsection{Problem instances} \label{sec:problem instances}

A problem instance is defined by the following quantities:
\begin{itemizedense}
    \item The set of parts $I$ with $|I|=\num{48}$. Each part $i \in I$ is assigned a value $V_i \in \mathbb{R}_{>0}$, which can be used to calculate the relative value $v_i := 100 V_i / V$ (in percentage) with the total value $V := \sum_{i \in I} V_i$. Furthermore, each part is assigned a volume $\mathcal{V}_i \in \mathbb{R}_{>0}$ (assuming a rectangular cuboid shape). The production dependencies between parts are defined by a set of tuples $\phi \in I^2$, where $(i,j) \in \phi$ indicates that the production of part $j$ (called ``parent'') requires part $i$ (called ``child''). Every child has only one parent, which leads to a dependency tree, the \gls{PBS}. The level $L_i := \sum_{(i',j) \in \phi} \delta_{i'i} (L_j+1) \in [0,4]$ of each part $i \in I$ represents the necessary number of production steps until the final product as assembled. 
    \item The set of production sites $K$ with $|K|=\num{43}$. Each site $k \in K$ is assigned one region $r \in R$ out of the set of all regions $V$ with $|V|=\num{18}$. To this end, we define $V_{vk} \in \{0,1\}$ with $\sum_{v \in V} V_{vk}=1$ for all $k \in K$, where $V_{vk} = 1$ indicates that region $v \in V$ is assigned to site $k \in K$ (and $V_{vk} = 0$ that it is not assigned). Furthermore, $K^{\text{min}}_{k} \in [0,100] \subset \mathbb{N}_0$ and $K^{\text{max}}_{k} \in [0,100] \subset \mathbb{N}_0$ define the minimum and maximum workshare (in percentage), respectively, for site $k \in K$. The production processes at the sites is not further concerned, it is, for example, presumed that there is always sufficient capacity for immediate production.
    \item The set of warehouses $W$ with $|W|=\num{28}$ that can be used as intermediate transportation hubs but don't produce anything.
    \item The set of suppliers $U$ with $|U|=\num{29}$. For each supplier $u \in U$, $U^{\text{min}}_{u} \in [0,100] \subset \mathbb{Z}$, $U^{\text{max}}_{u} \in [0,100] \subset \mathbb{N}_0$, and $U^{\text{target}}_{u} \in [0,100] \subset \mathbb{N}_0$ define the minimum, maximum, and target workshare (in percentage), respectively.
    \item The set of feasible site-supplier combinations $f_i \subseteq K \times U$ that are capable of producing part $i \in I$ with $|f_i| \geq 1$. The set of assignable sites for a part $i$ is denoted by $K_i := \{ k | k \in K \land \exists u \in U : (k,u) \in f_i \}$ with $|K_i| \geq 1$. Any part $i$ that can be assigned to at least two sites (\ie, $|K_i|\geq2$) can be used for double sourcing. Any other part can not be used for double sourcing. The set of assignable regions for a part $i$ is denoted by $V_i := \{ v | v \in V \land V_{vk} = 1 \land \exists u \in U : (k,u) \in f_i \}$ with $|V_i| \leq |K_i|$. Any part $i$ that can be assigned to at least two regions (\ie, $|V_i|\geq2$) can be double sourced in two distinct regions. Furthermore, for each tuple $(k,u) \in f_i$, the production time of the corresponding manufacturing process is given by $t_{i \rightarrow (k,u)} \in \mathbb{R}_{\geq0}$.
    \item The (possible empty) set of transportation methods $M^i_{p \Rightarrow q} \subset M_i \subseteq M$ for transporting part $i \in I$ from site or warehouse $p \in K \cup W$ to site or warehouse $q \in K \cup W$, where $|M|=\num{17073}$. The mappings $c^n := M \rightarrow \mathbb{R}_{\geq0}$ with $n \in \{1,2,3\}$ yield the carbon dioxide emissions ($c^1$), costs ($c^2$), and transportation time ($c^3$), respectively, for each transportation method. In addition, each transportation method is assigned a total cargo volume $\mathcal{V}_m \in \mathbb{R}_{>0}$ (assuming a rectangular cuboid shape). We denote the set of transportation methods for part $i \in I$ by $M_i$ and the set of all transportation methods by $M$.
    \item The primary source share $\alpha_i \in [0.5, 1]$ is the fraction of the production of part $i \in I$ at the primary production site. The secondary share is consequently given by $1-\alpha_i \in [0,0.5]$ and represents the fraction of the production of part $i$ at the secondary production site. Without loss of generality, the primary source share is by definition always larger than or equal to the secondary source share, \ie, $\alpha_i \geq 1-\alpha_i$ for all $i \in I$.
\end{itemizedense}

The data provided by Airbus contains a single logistics network for the production of an aircraft with the \gls{PBS} shown in \cref{fig:problemSketch}a and the production location and warehouses shown in \cref{fig:problemSketch}b. This single instance can be varied by a modification of the primary source shares (denoted $\alpha_i$). However, we mainly focus on the solution of a single instance version.

\subsection{Multi-Objective QCBO Model} \label{sec:qcbo model}

We use an index $n \in \{1,\dots,4\}$ to label the the four \glspl{KPI} as carbon dioxide emissions ($n=1$), costs ($n=2$) production time ($n=3$), and supplier target workshare fulfillment ($n=4$). We model these quantities in the following way:
\begin{enumeratedense}
    \item To determine the carbon dioxide emissions ($n=1$) and costs ($n=2$), we take all chosen transportation methods of all parts into account. For each mode of transportation, we count the emissions and pricing of the shipment in proportion to the relative shipping space occupied by the shipped part. This approach takes into account the fact that each part takes up only a fraction of the available cargo space, so that the remaining space can be filled with other goods and is not wasted. It has a significant impact on the cost of large-capacity transportation modes, such as ships, which would otherwise be counted with their full emissions even if only a tiny fraction of their cargo space is used to ship a part.
    \item For the production time ($n=3$), we consider the weighted cumulative transport times of all selected transportation methods instead of the actual lead time of the assembled product for two reasons. First, modeling the exact lead time greatly increases the complexity of the model. Second, we aim for a robust result, taking into account that in an efficient and resilient supply chain, all transportation modes should be chosen as small as possible independently of each other. As a weight of the transportation times we choose the level of the product in the \gls{PBS}. This choice recognizes that parts deeper in the tree have more (indirect) dependencies, so their transport time is more important and should result in a relatively larger contribution.
    \item For the supplier target workshare fulfillments ($n=4$), we choose the cumulative quadratic deviation from the relative target workshare for each supplier.
\end{enumeratedense}

Next, we describe how we realize this modeling approach for the \glspl{KPI}. To that end, we first discuss the shipment of parts. For each part $i \in I$, the available transportation methods $M_i \subset M$ and the joint set of production sites $K$ and warehouses $W$ can also be considered as a graph $G_i:=(K \cup W, M_i)$, where an edge $M^i_{p \Rightarrow q} \subset M_i \subseteq M$ is meant to understood to connect node $p$ to node $q$. In short, the graph $G_i$ defines how part $i$ can be moved between different locations.

Transportation methods can be combined into routes, where each route $r \in R^i_{k \Rightarrow l}$ from site $k \in K$ to site $l \in K$ represent a distinct path on the graph $G_i$, \ie, it consists of a sequence of $\nu \in \mathbb{N}$ transportation methods (graph edges) $r := (m_1, \dots, m_{\nu})$ with $m_1 \in M^i_{k \Rightarrow w}$, $m_s \in M^i_{w \Rightarrow w'}$ for $s \in \{2,\dots,\nu-1\}$, and $m_{\nu} \in M^i_{w \Rightarrow l}$ with $w, w' \in K \cup W$. For $\nu>1$, $r$ is an indirect route, where for $\nu=1$, $r := (m)$ with $m \in M^i_{k \Rightarrow l}$ is a direct route. By definition, routes always start and end at production sites, but they may involve transportation methods which pass through warehouses. Here, we use $R^i_{k \Rightarrow l}$ to denote the set of all available routes from $k$ to $l$ which visit each production site and warehouse (graph nodes) at most once along the path. Although this set can be quite extensive, it is not explicitly used in our final model. This will be clarified further below. If $R^i_{k \Rightarrow l} = \emptyset$, there is no available connection between $k$ and $l$. The set of all available routes for a part $i$ is denoted by $R^i := \cup_{k,l \in K} R^i_{k \Rightarrow l}$. A part $i$ is considered immobile, if $R^i = \emptyset$ (and mobile otherwise).

To take these different routing options into account in our modeling approach, we define a set of \emph{terminals} for each part at each site, which is to be understood as a set of labels that can be used to identify the different outgoing routes from this site to other sites, considering all possible direct and indirect routes. We define the set of terminals for part $i \in I$ at site $k \in K$ by $T^i_k:=\{1,\dots,t^i_k\}$ with the total number of terminals $t^i_k \geq 1$. First of all, we consider a mobile part (with $R^i \neq \emptyset$), for which the number of terminals $t^i_k$ is given by the maximum number of outgoing routes over all destinations such that $|T^i_k| = \max_{l \in K} | R^i_{k \Rightarrow l}|$. Then, there exists a surjective but not necessarily injective mapping $\tau^i_{kl}: T^i_k \rightarrow R^i_{k \Rightarrow l}$ for all $l \in K$ with $R^i_{k \Rightarrow l} \neq \emptyset$, such that each each terminal is assigned to a (not necessarily unique) route. In other words, given a part $i \in I$, a start site $k \in K$, and a destination site $l \in K$ with $R^i_{k \Rightarrow l} \neq \emptyset$, each terminal $t \in T^i_k$ is associated with a route $r \in R^i_{k \Rightarrow l}$. Conversely, all routes starting from site $k$ can be associated with at least one terminal in $t \in T^i_k$. We define a corresponding inverse mapping $\overline{\tau}^{i}_{kl} : R^i_{k\Rightarrow l} \rightarrow T^i_k$ from a route to a terminal, which is to be understood as an arbitrary but fixed realization of the inversion. We only presume that suitable mappings $\tau^i_{kl}$ and $\overline{\tau}^{i}_{kl}$ exist, but do not need to define them explicitly.

For an immobile part (\ie, $R^i = \emptyset$), there are by definition no outgoing routes, independent of where the part is produced. In this case, we presume that there is a single terminal that represents the absence of any transportation. Summarized, we define the the total number of terminals as
\begin{align} \label{eqn:numterminals}
    t^i_k := \begin{cases} \max_{i \in I,l \in K} | R^i_{k \Rightarrow l}| & , \text{if } R^i \neq  \emptyset \\ 1 & , \text{if } R^i = \emptyset \end{cases}
\end{align}
for any part $i \in I$, either mobile or immobile.

The concept of terminals allows us to integrate different routing options into the existing assignment task. For this purpose, in addition to the requirement that each part must be assigned to a site and supplier, we additional require that it must also be assigned to a terminal at the site. That means that a site-terminal-supplier tuple $(k,t,u)$ with $k \in K$, $t \in T^i_k$ and $u \in U$ must be assigned to each part $i \in I$ (both for the primary and secondary production). Based on the feasible site-supplier combinations $f_i$, the feasible site-terminal-supplier tuples are given by $F_i := \{ (k,t,u) \,|\, t \in T^i_k \,\forall\, (k,u) \in f_i \}$. 

Presume a site-terminal-supplier assignment $(k,t,u) \in F_i$ for part $i \in I$ and a site-terminal-supplier assignment $(l,s,v) \in F_j$ for part $j \in I$, where $(i,j) \in \phi$. Then, $\tau^i_{kl}(t) = r$ assigns a route $r:= (m_1,\dots,m_{\nu}) \in R^i_{k \Rightarrow l}$, which is used to transport part $i$ from site $k$ to site $l$, where part $j$ is produced. By definition, the assigned route $r$ depends on the sites $k$ and $l$ as well as the terminal $t$, but not on the suppliers $u$ and $v$, and also not on the terminal $v$. 

Each transportation method within the route $r$ implies a shipment that creates carbon dioxide emissions, increase costs and takes time. The total contribution of $r$ to the first three \glspl{KPI} is consequently given by the cumulative contributions of the transportation methods along the path, which we define as
\begin{align} \label{eqn:cover}
    c_{(k,t)\Rightarrow (l,s)}^{in} := \sum_{m \in \tau^i_{kl}(t)} c^n(m) \gamma^{in}(m)
\end{align}
for $R^i_{k \Rightarrow l} \neq \emptyset$, based on the the prescribed \gls{KPI} contributions $c^n$ from the problem instance and a custom scaling factor $\gamma^{in} \geq 0$. While $c_{(k,t)\Rightarrow (l,s)}^{in}$ is formally independent of the terminal $s$, we still include it in the notation for clarity. We use the scaling factor
\begin{align} \label{eqn:gamma}
    \gamma^{in}(m) := \begin{cases} \frac{\mathcal{V}_i}{\mathcal{V}_{m}} & , \text{if } n \in \{1,2\} \\ L_i & , \text{if } n=3\end{cases}
\end{align}
as explained above. Since every part has a non-vanishing volume and must fit into the cargo space, $\mathcal{V}_i / \mathcal{V}_{m} \in (0,1]$.

We aim to model the optimization problem as a \gls{QUBO}. To this end, we define a set of binary variables $y$, each identifiable by five indices. Specifically, the binary variable $y_{i\rightarrow(k,t,u)}^{a} \in \{0,1\}$ represents the assignment of part $i \in I$ to terminal $t \in T^i_k$ of site $k \in K$ and supplier $u \in U$, with $(k,t,u) \in F_i$. Here, $a \in A:=\{1,2\}$ indicates either whether the assignment refers to the the primary production ($a=1$) or the secondary production ($a=2$). If a part $i \in I$ does not allow double sourcing (\ie, $|K_i|=1$), we set $y_{i\rightarrow(k,t,u)}^{a=2} := y_{i\rightarrow(k,t,u)}^{a=1}$, effectively reducing the decision space. 

With these premises, we can formulate the optimization problem as:

\begin{condstrip}
\begin{subequations}
\begin{align}
\min_{y} & (C_1(y), C_2(y), C_3(y), C_4(y)) \label{eqn:model1} \tag{\theparentequation} \\
    \text{s.t. } & \sum_{a,b \in A, u, v \in U} \mu_i^{ab} y_{i\rightarrow(k,t,u)}^{a} y_{j\rightarrow(l,s,v)}^{b} = 0 \,\forall\, (i,j) \in \phi \,\forall\, k,l \in K \,\forall\, t \in T^i_k \,\forall\, s \in T^j_l \land R^i_{k \Rightarrow l} = \emptyset \label{eqn:model1:cI} \\
     & \sum_{(k,t,u) \in F_i} y_{i\rightarrow(k,t,u)}^{a} = 1 \,\forall\, i \in I \,\forall\, a \in A \label{eqn:model1:cII} \\
     & \Big( \sum_{(k',t,u) \in F_i} \delta_{k'k} y_{i\rightarrow(k,t,u)}^{a=1} \Big) \Big( \sum_{(k',t,u) \in F_i} \delta_{k'k} y_{i\rightarrow(k,t,u)}^{a=2} \Big) = 0 \,\forall\, i \in I \,\forall\, k \in K \land |K_i|\geq2 \label{eqn:model1:cIII} \\
     & \Big( \sum_{(k,t,u) \in F_i} V_{vk} y_{i\rightarrow(k,t,u)}^{a=1} \Big) \Big( \sum_{(k,t,u) \in F_i} V_{vk} y_{i\rightarrow(k,t,u)}^{a=2} \Big) = 0 \,\forall\, i \in I \,\forall\, v \in V \land |V_i|\geq2 \label{eqn:model1:cIV} \\
     & \sum_{i \in I, (k',t,u) \in F_i, a \in A} \delta_{k'k} v_i \alpha_i^a y_{i\rightarrow(k,t,u)}^{a} \in [K^{\text{min}}_{k},K^{\text{max}}_{k}] \,\forall\, k \in K \label{eqn:model1:cV} \\
     & \sum_{i \in I, (k,t,u') \in F_i, a \in A} \delta_{u'u} v_i \alpha_i^a y_{i\rightarrow(k,t,u)}^{a} \in [U^{\text{min}}_{u},U^{\text{max}}_{u}] \,\forall\, u \in U \label{eqn:model1:cVI}
\end{align}
\end{subequations}
\end{condstrip}

In summary, this is a multi-objective \gls{QCBO} with four \glspl{KPI} and six types of constraints. As \glspl{KPI}, we use carbon dioxide emissions ($C_1$), costs ($C_2$), and times ($C_3$) of all used means of transportation for all parts, with the former considered in proportion to the relative volume of the load, and the latter considered in proportion to the level of each part as a measure of its priority. The fourth \gls{KPI} is the supplier workshare target fulfillment ($C_4$), taken into account as the quadratic deviation from the target workshare of each supplier. Formally, the \glspl{KPI} are defined by

\begin{subequations} \label{eqn:C1234b}
\begin{align} \label{eqn:C123b}
    C_n(y) := \frac{1}{d_n} \sum_{\substack{(i,j) \in \phi \\(k,t,u) \in F_i,a \in A \\ (l,s,v) \in F_j, b \in A \\ R^i_{k \Rightarrow l} \neq 0}} & \alpha_i^a c_{(k,t)\Rightarrow(l,s)}^{in} \\
    & \times y_{i\rightarrow(k,t.u)}^{a} y_{j\rightarrow(l,s,v)}^{b}  \nonumber 
\end{align}
\end{subequations}

for dioxide emissions ($n=1$), costs ($n=2$), and transportation times ($n=3$), whereas the supplier target workshare fulfillments is defined by

\begin{subequations}
\begin{align} \label{eqn:C4b}
    C_4(y) := \frac{1}{d_4} \sum_{u \in U} \Big[ & \sum_{\substack{i \in I \\ (k,t,v) \in F_i, a \in A}} \!\!\!\!\!\! \delta_{uv} v_i \alpha_i^a y_{i\rightarrow(k,t,u)}^{a} \\
    & - U^{\text{target}}_{u} \Big]^2, \nonumber 
\end{align}   
\end{subequations}
where we recall the relative value $v_i$. Furthermore, we make use of the abbreviation
\begin{align}
    \alpha_i^a := \begin{cases} \alpha_i & , \text{if } a=1 \\ 1-\alpha_i & , \text{if } a=2\end{cases}
\end{align}
for the primary (secondary) source share $\alpha_i$ ($1-\alpha_i$) and introduce the empirically chosen rescaling factor
\begin{align}
    d_n := \begin{cases} \hat{c}^n & ,\text{if } n \in \{1,2\} \\ \hat{c}^{n=3} \hat{L} & ,\text{if } n = 3 \\ 100 & , \text{if } n = 4\end{cases}
\end{align}
to align the magnitudes of the \glspl{KPI} without altering the model structure. It contains the maximum contributions
\begin{align}
    \hat{c}^n := \max_{\substack{i \in I \\ p,q \in K \cup W \\ m \in M^i_{p \Rightarrow q}}} c^n(m)
\end{align}
and the level mid-range
\begin{align}
    \hat{L} := \frac{\min_{i \in I} L_i + \max_{i \in I} L_i}{2} = 2,
\end{align}
which is used as an additional factor for the transportation times based on the observation that the scaling factor from \cref{eqn:gamma} does not contain a relative value for $n=3$. The constant value \num{100} for $n=4$ is motivated by the percentage measure of the target workshares. By construction, all \glspl{KPI} are non-negative.

The six types of constraints ensure the requirements of a feasible supply chain: avoid non-existing connections (\cref{eqn:model1:cI}), allow one one assignment for each production source (\cref{eqn:model1:cII}), enforce different sites for primary and secondary source (\cref{eqn:model1:cIII}), enforce different regions for primary and secondary sources (\cref{eqn:model1:cIV}), respect the site workshare limits (\cref{eqn:model1:cV}), and respect the supplier workshare limits (\cref{eqn:model1:cVI}).

\subsection{QUBO Model} \label{sec:qubo model}

To transform \cref{eqn:model1} into a \gls{QUBO}, we perform the following steps:
\begin{enumeratedense}
    \item \Cref{eqn:model1} involves multiple objectives. To reduce the multi-objective assignment problem to a single-objective assignment problem, we perform a scalarization of the \glspl{KPI} with the weights $w := (w_1, w_2, w_3, w_4) \in [0,1]^4$, where $\sum_n w_n = 1$, such that the new optimization objective reads $\sum_{n=1}^4 w_n C_n(y)$. Different choices of weights $w$ represent different compromises between the objectives in the usual sense.
    \item The explicit use of terminals results in a large number of terms in some expressions, which make them computationally challenging to handle. To address this, we employ a two-step preprocessing of the instance: first, a pathfinding step using Dijkstra's algorithm to eliminate sub-optimal terminals (based on the chosen scalarization weights $w$), and second, a refinement of the feasible solution space to reduce the overall number of variables (indepentent of the weights $w$). The preprocessing is explained in detail in \cref{app:preprocessing}. As a result, we obtain a reduced model with less terms and variables.
    \item The reduced model still contains constraints, which need to be transformed into penalty terms. For the equality constraints, this can be realized straightforwardly. The box constraints, however, must first be converted into equality constraints. To achieve this, we employ a standard method~\cite{babbush2013}, which requires us to approximate the relatives values $v_i$ and the primary source shares $\alpha_i$ as rational numbers for all parts $i\in I$. Using auxiliary (ancilla) variables $z \in \{0,1\}^{N_z}$, we can then reformulate the box constraints as equality constraints and transform them into penalty terms. This method is described in \cref{app:box constraints}. 
\end{enumeratedense}

As a result, we arrive at the final \gls{QUBO} form
\begin{subequations} \label{eqn:qubo}
\begin{align}
    \min_{x} Q(x)
\end{align}
with the objective
\begin{align} \label{eqn:qubo:Q}
    Q(x) & := Q(w,\lambda,x)\\
    & := \sum_{n=1}^4 \overline{C}_n(w,\overline{y}) + \sum_{n=1}^6 \lambda_n P_n(x) \nonumber ,
\end{align}
\end{subequations}
which consists of the four weighted \glspl{KPI} $\overline{C}_n(w_n,\overline{y})$ and the six penalty terms $P_n(x)$ with the corresponding Lagrange multipliers $\lambda_n > 0$, where we use the notation $\lambda := (\lambda_1,\dots,\lambda_6)$. The set of binary variables $x \in \{0,1\}^{N_x}$ with $N_x := N_y + N_z$ fully describes the supply chain configuration. All expressions from \cref{eqn:qubo:Q} are defined in \cref{app:qubo}.

The resulting \gls{QUBO} objective, \cref{eqn:qubo:Q}, can also be written in a quadratic form
\begin{align} \label{eqn:quadratic}
    Q(x) = x^{\intercal} Q x
\end{align}
with a matrix $Q:=Q(w,\lambda)\in\mathbb{R}^{N_x \times N_x}$. Alternatively, using the coordinate transformation $s := 2x-1$, the \gls{QUBO} can also be straightforwardly converted into an Ising optimization problem~\cite{lucas2014}
\begin{subequations} \label{eqn:ising}
\begin{align} \label{eqn:ising:opt}
    \min_{s} H(w,\lambda,s)
\end{align}
based on the Ising model (also called Ising Hamiltonian)
\begin{align} \label{eqn:ising:H}
    H(w,\lambda,s) := H(s) = s^{\intercal} J s + h s
\end{align}
\end{subequations}
with spin variables $s \in \{-1,1\}^{N_x}$, a matrix $J:=J(w,\lambda)\in\mathbb{R}^{N_x \times N_x}$, and a vector $h:=h(w,\lambda)\in\mathbb{R}^{N_x}$.

\section{Methods} \label{sec:methods}

\begin{figure*}[!tb]
	\centering
	\includegraphics[]{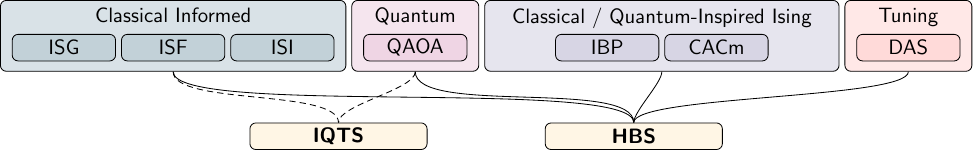}
	\caption{Sketch of the algorithmic architecture. We propose two hybrid quantum-classical solvers, \acrshort{IQTS} (\cref{sec:iqts}) and \acrshort{HBS} (\cref{sec:hbs}), which make use of various components: \acrshort{ISG} (\cref{sec:isg}), \acrshort{ISF} (\cref{sec:isf}), \acrshort{ISI} (\cref{sec:isi}), \acrshort{QAOA} (\cref{sec:qaoa}), \acrshort{IBP} (\cref{sec:ibp}), \acrshort{CACm} (\cref{sec:cacm}), \acrshort{DAS} (\cref{sec:das}).}
	\label{fig:solvercomponents}
\end{figure*}

To solve \cref{eqn:qubo}, we propose two hybrid quantum-classical solvers: \acrfull{IQTS} and \acrfull{HBS}. While \Gls{IQTS} leverages the domain knowledge of the problem structure by integrating \gls{QAOA} within an exploratory tree decomposition, \gls{HBS} is a bilevel optimization framework with a highly modular architecture, designed to combine various methods such as \gls{QAOA}, \gls{IBP}, \gls{CACm}, and \gls{DAS} to enhance flexibility, performance, and scalability. Both methods rely on similar solver components, as sketched in \cref{fig:solvercomponents}. In the following, we first describe these components, which then allows us to present \gls{IQTS} and \gls{HBS}.

\subsection{Solver Components} \label{sec:components}

In this section, we provide a brief summary over all seven solver components, which can be grouped into four categories:
\begin{itemizedense}
    \item Three classical procedures based on our structural knowledge about the problem, which is why we also call them ``informed'' procedures: \Gls{ISG}, \gls{ISF}, and \gls{ISI}.
    \item A quantum optimization algorithm, \gls{QAOA}, which is fine-tuned for our problem.
    \item Two classical Ising solvers, \gls{IBP} and \gls{CACm}. The latter can also be considered as a quantum-inspired approach.
    \item A derivative-free hyperparameter tuning procedure, \gls{DAS}.
\end{itemizedense}
All of these components are heuristic algorithms.

\subsubsection[ISG]{Informed Solution Generator (ISG)} \label{sec:isg}

Since we have a detailed knowledge about the problem structure, we can use this knowledge to generate randomized feasible solutions (that fulfill all constraints). For that, we propose the following procedure, which we call \acrfull{ISG}:
\begin{enumeratedense}
    \item Start with all parts unassigned.
    \item Iteratively draw an unassigned part $i \in I$ from the lowest available \gls{PBS} level (the reason for this order is that lower levels tend to have fewer available options).
    \item Collect all feasible site-supplier pairs $g_i$ for the part $i$ that satisfy the constraints, given the assignments made previously. Randomly assign the part $i$ to one of these options. Repeat this step for the primary and secondary sources.
    \item Repeat from step 2 until all parts are assigned, resulting in $\overline{y}$. Calculate the correct values for $z$ and return the corresponding solution $x = (\overline{y},z)$, which is a feasible by construction.
\end{enumeratedense}

\subsubsection[ISF]{Informed Solution Fixer (ISF)} \label{sec:isf}

We can also exploit the problem structure to project an infeasible solution back to the feasible domain without disturbing it more than necessary. Our proposed procedure \acrfull{ISF} works as follows:
\begin{enumeratedense}
    \item Start with $x$, a candidate solution to \cref{eqn:qubo} that violates the constraints.
    \item Iterate through all sites $k \in K$ and suppliers $u \in U$ and mark whether the workshare is overassigned (too much workshare) or underassigned (not enough workshare). If overassigned, randomly unassign parts until the workshare is satisfied. This step can lead to an update of $\overline{y}$.
    \item Iterate through all parts $i \in I$ of decreasing \gls{PBS} level and mark whether they are unassigned or wrongfully assigned (\ie, in a way that violates the constraints). Try to reassign each wrongfully assigned part: If there are underassigned sites/suppliers, try to reassign the mismatched products to those sites/suppliers and avoid overassigned sites/suppliers. Reassignments must meet all constraints. If a reassignment cannot be resolved, skip it. Again, this step can lead to an update of $\overline{y}$.
    \item Repeat from step 2 until all assignments are feasible, resulting in $\overline{y}$, or until a certain number of iterations have been performed, at which point the method has ran out of budget and no fixed solution can be returned. In case that a set of feasible assignments has been found, calculate the correct values for $z$ and return a fixed solution $x = (\overline{y},z)$, which is feasible by construction.
\end{enumeratedense}

\subsubsection[ISI]{Informed Solution Improver (ISI)} \label{sec:isi}

Finally, we can also use our knowledge to iteratively improve the objective of a feasible solution without leaving the feasible domain. For that, we propose a procedure called \acrfull{ISI}, which works in the following way:
\begin{enumeratedense}
    \item Start with $x$, a feasible candidate solution to \cref{eqn:qubo}.
    \item Randomly draw a part $i \in I$.
    \item Find all feasible combinations of primary source and secondary source assignments for part $i$.
    \item Iterate through the feasible options in a random order and, if they would lead to an improvement in the optimization objective, reassign the part $i$ accordingly, which leads to an update of $\overline{y}$. In case of a reassignment, stop exploring the feasible options with a predefined probability or continue otherwise. Stop when all feasible options have been explored.
    \item Repeat from step 2 until a certain number of iterations has been reached. Calculate the correct values for $z$ and return an improved $x = (\overline{y},z)$, which is feasible by construction and leads to an objective at least as good as the original candidate solution.
\end{enumeratedense}

\subsubsection[QAOA]{Quantum Approximate Optimization Algorithm (QAOA)} \label{sec:qaoa}

\Gls{QAOA}~\cite{farhi2014,blekos2024} is a hybrid quantum-classical heuristic to solve Ising optimization problems of the form of \cref{eqn:ising} on a gate-based quantum computer. As a special form of \gls{VQA}, the key idea is to optimize the parameters of a parameterized quantum circuit iteratively using a classical optimizer until convergence is achieved. The measurement results of the quantum device then correspond to candidate solutions $s$ of \cref{eqn:ising:opt}.

We make use of a variant of \gls{QAOA}~\cite{farhi2014} with a problem-specific setup:
\begin{itemizedense}
    \item Mixer: We choose a $XY$-mixer~\cite{wang2020}, which consists of correlated Pauli-x and Pauli-y terms $\frac{1}{2} ( \hat{\sigma}^i_x \hat{\sigma}^j_x + \hat{\sigma}^i_y \hat{\sigma}^j_y )$, for every combination of variables $(i,j)$ that correspond the an assignment of the same product for the same source. This mixer enforces the constraint $P_2(x)$. In addition, we choose the uncorrelated Pauli-x terms$\hat{\sigma}^i_x$ for all ancilla variables.   
    \item Initialization: The initial state is chosen to be a W state~\cite{cruz2019,mukherjee2020} for every set of variables that correspond the an assignment of the same product for the same source in analogy to the mixer choice. W states can be prepared with logarithmic complexity~\cite{cruz2019}. The state for all ancilla variables is prepared in the usual way with a Hadamard gate.
    \item Schedule: Motivated by recent findings~\cite{montanezbarrera2024}, we make use of a fixed linear ramp \gls{QAOA} protocol, for which the \gls{QAOA} parameters~\cite{farhi2014}, typically called $\gamma$ and $\beta$, are set to fixed values. This means that no additional classical optimization loop is necessary to tune these parameters, which greatly reduces the required quantum hardware access time. Specifically, for a depth $p \geq 2$, we choose constant cost parameters $\gamma_i =  (i-1)/(p-1)$ and constant mixer parameters $\beta_i = 1-\gamma_i$ for $i \in [2,p]$. For the special case $p=1$, we define $\gamma_1=\beta_1=\frac{1}{2}$.
\end{itemizedense}
This setup requires $n$ qubits to solve a \gls{QUBO} with $n$ binary variables. The only hyperparameter is the depth $p \in \mathbb{N}$. As the size of our problem formulation, \cref{eqn:qubo}, is too large to be solved with the quantum hardware available to us, we apply \gls{QAOA} only on sub-problems. While \gls{QAOA} is genuinely a hybrid quantum-classical algorithm, the choice of fixed parameters reduce our variant to a pure quantum algorithm.

\subsubsection[IBP]{Iterative Belief Propagation (IBP)} \label{sec:ibp}

\Gls{BP} is an algorithm which approximately samples from a high dimensional discrete probability distribution when the probability function can be expressed as a product of factors that involve a small subset of the variables~\cite{Braunstein2005_SP,Pearl1982_BP}. If the graph that connects factors and variables (called the factor graph) forms a tree, \gls{BP} is know to converge exactly in polynomial time. In the context of \glspl{QUBO}, the probability distribution we would like to sample from is the Boltzmann distribution for a given inverse temperature $\beta$, defined as $P(x) \propto e^{-\beta Q(x)}$ where $Q(x)$ is the \gls{QUBO} cost function of the form of \cref{eqn:qubo:Q}. Because $Q(x)$ is a linear combination of terms involving edges in the \gls{QUBO} graph these edges will then become factors in the graph used by \gls{BP}. Thus, if we have a \gls{QUBO} with tree connectivity, \gls{BP} provides a polynomial time algorithm for finding the ground state (by setting $\beta$ to a very large value). 

\Gls{BP} can be outlined as follows. Consider a \gls{QUBO} in quadratic form, \cref{eqn:quadratic}, with $N_x$ variables $x$ and a corresponding matrix $Q$. For each ordered pair of variables $i,j \in \{1,\dots,N_x\}$ such that $Q_{ij} \neq 0$,  \gls{BP} stores two ``messages'' $\mu^{a}_{ij} \in [0,1]$ with $a \in \{0,1\}$. These messages are typically initialized randomly and updated iteratively according to the classical belief propagation equations~\cite{mezard2009information}:

\begin{subequations}
\begin{equation}\label{eqn:bp2}
    \mu^{a}_{ij} \rightarrow \frac{(\mu^{a}_{ij})^{*}}{\sum_{b \in \{0,1\}} (\mu^{b}_{ij})^{*}} 
\end{equation}
with
\begin{equation}\label{eqn:bp1}
    (\mu^{a}_{ij})^{*} := \!\!\prod_{\substack{k \in \mathcal{N}(i)\\k \neq j}} \!\left[\sum_{b \in \{0,1\}}e^{-\beta Q_{ik} ab} \!\!\prod_{\substack{l \in \mathcal{N}(k)\\ l \neq j}} \!\! \mu^{b}_{lk} \right].
\end{equation}
\end{subequations}
The set $\mathcal{N}(i)$ refers to all of the connected QUBO variables (which can include $i$ itself). Once the messages have converged, they can be used to compute marginal probabilities which allow us to sample from the desired (Boltzmann) distribution.

On the other hand, \gls{SA} provides a different approach for sampling from the Boltzmann distribution of a \gls{QUBO}~\cite{Kirkpatrick1983_SA}. In \gls{SA}, the current solution $x$ is iteratively updated by choosing a random variable and then flipping it according to the Metropolis Hastings criterion. This can also be interpreted as choosing a random sub-problem of the full problem (of just one variable in this case) and sampling from the Boltzmann distribution of that sub-problem while keeping all other variables fixed. Unlike \gls{BP}, \gls{SA} has the useful property that it works for any \gls{QUBO} connectivity given enough iterations, however in some cases convergence time can be exponential in the problem size.

The proposed \gls{IBP} procedure attempts to bridge the gap between \gls{BP} and \gls{SA} as follows~\cite{Reifenstein2024BP}. Given a \gls{QUBO} graph with sparse connectivity, it is possible to choose a subset of variables such that the subgraph formed by them is a tree. If the given \gls{QUBO} graph is sparse and more tree-like, these sub-graphs will be larger, whereas if the graph has full connectivity these subgraphs will have size 2 at most. \Gls{IBP} works by randomly choosing one of these subgraphs (which can be done efficiently) and applying \gls{BP} on it. This allows us to sample from the Boltzmann distribution for a subset of variables similar to \gls{SA}. However, unlike \gls{SA}, many variables will be updated at once. Iteratively choosing these sub-trees and updating $x$ using \gls{BP} will allow us to construct a Markov chain which converges to the Boltzmann distribution as desired. A convenient property of \gls{IBP} is that in the limiting case of a \gls{QUBO} that is already a tree, the algorithm will be identical to regular \gls{BP}, and thus converge in polynomial time (given the hyper-parameter $\beta$ is chosen correctly). On the other hand, if the problem is fully connected, \gls{IBP} will reduce to a slightly modified version of \gls{SA} in which two variables are updated at once instead of one. Similar to \gls{SA}, \gls{IBP} is parameterized by the hyper-parameter $\beta$ which can be annealed for better performance. Many problems of industrial and academic interest lie somewhere in-between these two extremes thus it is reasonable to postulate that \gls{IBP} can outperform both \gls{BP} and \gls{SA} in some practically relevant cases. In particular, since the optimization problem studied in this work has an inherent tree structure to it from the \gls{PBS}, it is reasonable that \gls{IBP} may provide some improvement in solution speed and accuracy by exploiting this structure.

\subsubsection[CACm]{Chaotic Amplitude Control with Momentum (CACm)} \label{sec:cacm}

The \gls{CACm} algorithm belongs to the class of Ising solvers based on \glspl{ODE}~\cite{leleu2019destabilization,leleu2021scaling}. It has been shown to perform competitively against other state-of-the-art solvers such as \gls{SA} and \gls{PT}~\cite{leleu2021scaling,Reifenstein2024,Leleu2024}, but its performance strongly depends on the choice of its hyperparameters~\cite{Reifenstein2024}. Presuming an Ising optimization problems of the form of \cref{eqn:ising}, \gls{CACm} is based on the concept of relaxing the binary space of $s \in \{-1,1\}^{N_x}$ to the real space $m \in \mathbb{R}^{N_x}$, which leads to a relaxed problem of the form $\min_m H(m)$.

The operator $F^1$, which represents the flow map of the deterministic trajectory of \gls{CACm}, is defined as follows~\cite{Leleu2024}:
\begin{subequations}
\begin{align}
\gamma \frac{d^2 u}{dt^2} + \frac{d u}{dt} = - \lambda u - \beta e \circ \nabla_m V,
\end{align}
\begin{align}
\frac{d e}{dt} = -\xi (m \circ m- a) \circ e,
\end{align}
\end{subequations}
where the symbol $\circ$ denotes the Hadamard product, $\nabla_m V$ the gradient of $V$ with respect to the vector $m$, $\beta$ a positive parameter and $e$ a vector of positive auxiliary variables. The variable $u$ is a real-valued vector representing the internal state of the \gls{CACm} and $t$ the continuous time. The term proportional to $\gamma$ represent the momentum which has been utilized for non-convex optimization~\cite{kalinin2023analog}. In practice, the operator $F^1$ used the bilevel optimization scheme for the contribution of \gls{CACm} represents the map given as $F^1: s_n : = \sign[m(t=0)] \rightarrow s_{n+1} : = \sign[m(t=T)]$.

In the current work, the potential $V$ is equal to the Ising Hamiltonian of the model $H$, \cref{eqn:ising:H}, relaxed to the continuous domain. The vector $\theta := (\lambda_1,\lambda_2,\gamma,\beta,\xi,a,T)$ contains the hyperparameters, which are briefly explained in the following. First of all, $\lambda_1$ and $\lambda_2$ are the initial and final values of the parameter $\lambda$ that is taken to be a linear function of time to perform ``annealing.''.  Furthermore, $\gamma$ represents the momentum and $\beta$ a positive parameter. We set $m_i := \phi(\tilde{\beta} u_i), \forall i \in \{1, \cdots, N_x\}$, where $\phi(m) := \frac{2}{1+e^{-m}} - 1$ is a sigmoidal function normalized to the domain $m \in [-1,1]$. Moreover, $\xi$ represents the speed of the error correction dynamics and $a$ represents target amplitude of the $m$ variables. Finally, $T$ denotes the number of steps of a single run of \gls{CACm}.

\subsubsection[DAS]{Dynamic Anisotropic Smoothing (DAS)} \label{sec:das}

\Gls{DAS} is a recently proposed derivative-free optimization approach for automated hyperparameter tuning. Presume that an algorithm has a hyperparameter space $\theta$. \Gls{DAS} involves iteratively sampling a set of solutions $\mathcal{L}[\theta]:= \{\theta_1, \theta_2, \cdots, \theta_R\}$, using the Ising solvers to infer an optimal value of $\theta$. This is achieved by iterating the following equations~\cite{Reifenstein2024}:
\begin{subequations} \label{eq: SDE}
\begin{align}
    \frac{dL}{dt} &= \alpha_L \left(LL^{\top} \frac{\partial h(L,\theta)}{\partial L} + \lambda L + \eta_L \right), \label{eq:dL} \\
    \frac{d\theta}{dt} &= \alpha_\theta \left(LL^{\top} \frac{\partial h(L,\theta)}{\partial \theta} + \eta_\theta\right), \label{eq:dtheta}
\end{align}
\end{subequations}
The variables $\theta$ and $L$ represent a distribution in the hyper-parameter space in which the trajectories are sampled from, where the vector $\theta$ and matrix $L$ correspond to the center and covariance of the distribution respectively. $R$ is the number of hyperparameters sampled per step of \Gls{DAS}. Equations~(\ref{eq:dtheta}) and~(\ref{eq:dtheta}) describe the gradient descent dynamics of the optimal parameter estimate $\theta$, as well as the estimation of the local curvature $L$ of the algorithm's performance in the hyperparameter space at $\theta$. The function $h$ and the hyperparameter $\lambda$ are described in more detail in~\cite{Reifenstein2024}.

In each step $n$ of \gls{DAS}, an ensemble of parameters $\mathcal{L}[\theta^{\gamma}]$ is updated conditionally on set of solution $\mathcal{S}^{\gamma}_{n}$ found by the Ising solver as follows:
\begin{align} \label{eq:BOA2}
\mathcal{L}[\theta_{n+1}^{\gamma}] = G(\mathcal{L[\theta}_n^{\gamma}]; \mathcal{S}^{\gamma}_{n}).
\end{align}
\noindent In \gls{DAS}, the operator $G$ estimates the gradient and local curvature of $\mathcal{C}_n$ (that is, \gls{DAS} is an approximate second order method).

\subsection[IQTS]{Informed Quantum-Enhanced Tree Solver (IQTS)} \label{sec:iqts}

Our first hybrid quantum-classical solver focuses on knowledge about the problem structure and how it can be used to decompose the problem into smaller sub-problems, which can then be tackled using quantum optimization. Since the decomposition approach is centered around the tree structure of the \gls{PBS}, we call this approach \acrfull{IQTS}. The method is largely I/O-bound and can achieve a very fast convergence by keeping the exploration close to the feasible domain.

The task of the solver is to provide a heuristic solution to the proposed \gls{QUBO}, \cref{eqn:qubo}. Given a user-defined subtree size $m$ and variable count $n$, the solver workflow can be summarized in the following way, where we recall the solver components from \cref{sec:components}:
\begin{enumeratedense}
    \item Start with a candidate solution $x$ from \gls{ISG}.
    \item Randomly select a part $i \in I$ from the \gls{PBS} that has not yet been selected, always prefer the part with the highest available level.
    \item For each part $i$, perform the following steps:
    \begin{enumeratedense}
        \item Randomly select a sub-tree of size $m$ from the \gls{PBS} that contains part $i$.
        \item Randomly select $n$ variables for the assignment of parts, which are contained in the sub-tree. Set all remaining variables to zero in order to nullify the product assignments.
        \item Build a sub-problem from the selected variables and solve it with \gls{QAOA}. Due to the selection of the sub-problem based on a sub-tree of the \gls{PBS}, the sub-problem variables are highly correlated.
        \item If necessary, repair $x$ using \gls{ISF}. This step can lead to an update of $x$.
        \item Perform a predefined number of \gls{ISI} iterations. This step can lead to an update of $x$.
    \end{enumeratedense}
    \item If all parts have been randomly selected once in step 2, make all of them freely available again.
    \item Repeat from step 2 until $\kappa$ repetitions have been reached, then return $x$.
\end{enumeratedense}

\subsection[HBS]{Hybrid Quantum-Classical Bilevel Solver (HBS)} \label{sec:hbs}

The second solver we propose is a high performance framework that leverages multiple sub-solvers in a collaborative way, which allows us to make use of a combination of classical, quantum-inspired and pure quantum optimizers. The framework is mostly CPU-bound and focused on computational efficiency, scalability, and generalizability for tackling complex models with bilevel optimization. We call this approach \acrfull{HBS}.

Specifically, the task of \gls{HBS} is to provide a heuristic solution to an Ising optimization problem of the form of \cref{eqn:ising}. We outline the framework in the following, where we recall the components from \cref{sec:components}. First of all, \gls{HBS} makes use of three specialized Ising solvers:
\begin{itemizedense}
    \item \Gls{CACm}~\cite{Leleu2024,leleu2019destabilization}, as a general-purpose, state-of-the-art quantum-inspired method,
    \item \gls{IBP}~\cite{Reifenstein2024BP}, as an iterative message-passing technique, designed to quickly approximate solutions for problems with a tree-like structure, and
    \item \gls{QAOA}~\cite{farhi2014quantum,mezard2009information} to include the exploitation of quantum effects.
\end{itemizedense}
Each solver leverages different properties of the problem instances, and together, they maximize the advantages of each approach. We denote the solvers by $\mathsf{Alg}^{\gamma}$ with $\gamma \in \{\text{\acrshort{CACm}, \acrshort{IBP}, \acrshort{QAOA}\}}$. In addition, \gls{HBS} also uses \gls{ISG} to generate solution candidates, \gls{ISF} to improve the feasibility of solution candidates and \gls{DAS} to tune hyperparameters. \Cref{{alg:hbs}} provides a summary of the framework.

\begin{algorithm}[!htb]
\caption{\Acrlong{HBS}} \label{alg:hbs}
\begin{algorithmic}[1]
    \State \textbf{Input:} Ising problem, \cref{eqn:ising}, with Ising Hamiltonian $H(s)$
    \State Initialize $s_0 \in \mathcal{S}_{0}$ using \acrshort{ISG}
    \While{not converged}
    \For{$\gamma \in \{\text{\acrshort{CACm}, \acrshort{IBP}, \acrshort{QAOA}\}}$}
    \State Refine solution: $s_{n+1} = F^{\gamma}(s_n, \theta_n^{\gamma})$
    \State Generate candidate set: $\mathcal{S}_{n+1}^{\gamma}$
    \EndFor
    \State Combine solutions: $\mathcal{S}_{n+1} = \mathbf{S}\left[\bigcup_{\gamma} \mathcal{S}_{n+1}^{\gamma}\right]$
    \State Tune: $\theta_{n+1}^{*} \approx \arg \min_{\theta} u(\mathcal{H}^*-\mathcal{H}(s_{n}))$ via \acrshort{DAS} step
    \State given as $\mathcal{L}_{n+1}[\theta] = G(\mathcal{L}_n[\theta]; \mathcal{S}_{n})$
    \State Fix solutions $s_{n+1} \in \mathcal{S}_{n+1}$ using \acrshort{ISF}
    \State \textbf{Record} candidate solutions $\mathcal{S}_{n+1}$
    \EndWhile
\end{algorithmic}
\end{algorithm}

Initially, \gls{ISG} is used to generate a candidate solution $s_0 \in \mathcal{S}_{0} \subset \{-1,1\}^{N_x}$, where $\mathcal{S}_{0}$ denotes the feasible set of states generated by \gls{ISG}. Then, a sequence of iterations $n \geq 1$ follows, during which the solution candidate is progressively refined until convergence is reached—that is, when further iterations no longer improve the solution quality. The task of each algorithm $\mathsf{Alg}_{\gamma}$ is to refine a solution from an initial guess $s_n \in \mathcal{S}_n$ to an improved solution $s_{n+1} \in \mathcal{S}^{\gamma}_{n+1}$, ensuring $H(s_{n+1}) \leq H(s_{n})$. Formally, this action is represented as a operator $F^{\gamma}$ given as follows:
\begin{align} \label{eq:BOA1}
    s_{n+1} := F^{\gamma}(s_n, \theta_n^{\gamma}), \forall s_n \in S_n, \forall \gamma,
\end{align}
where $\theta_n^{\gamma}$ denotes the hyperparameters of solver $\gamma$ at step $n$. In the case of \gls{CACm} for example, the operator $F^{\mathrm{CACm}}$ consists in the flow map of a deterministic trajectory represented by an \gls{ODE} as described in \cref{sec:cacm}.

The solutions found by all algorithms are then recombined to be the starting point for the next iterations, which involves sampling the subset of best solutions found among the union set of all solutions, as follows:
\begin{align}
    \mathcal{S}_{n} := \mathbf{S}\left[\bigcup_{\gamma} \mathcal{S}^{\gamma}_{n}\right],
\end{align}

\noindent where $\mathbf{S}$ is a sampling operation selecting with uniform distribution from candidate solutions with smallest Ising Hamiltonian $H(s_{n+1})$.

\begin{figure}[!htb]
	\centering
	\includegraphics[]{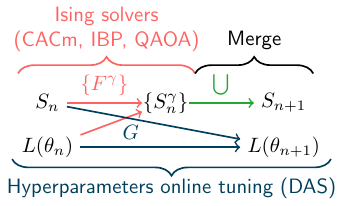}
	\caption{Sketch of how \gls{HBS} combines Ising solvers (\gls{CACm}, \gls{IBP}, and \gls{QAOA}) with hyperparameter tuning (\gls{DAS}).} \label{fig:BOA}
\end{figure}

We then adopt a bilevel optimization approach, where the parameters $\theta^{\gamma}_n$ are adaptively tuned to maximize the quality of the solution found. Specifically, we seek
\begin{align}
    \theta^{\gamma*}_n := \arg \min_{\theta} \mathcal{C}^{\gamma}_n(\theta),
\end{align}
where 
\begin{align}
    \mathcal{C}^{\gamma}_n(\theta) := u(H^*-H(F^{\gamma}(s_{n-1}, \theta)))
\end{align}
is a cost associated with the solution quality found after one step of the Ising solver and based on the the Heaviside step function $u$. The reference energy $H^*$ is defined as the minimum objective of \cref{eqn:ising} found so far. Since the functional dependency of $\mathcal{C}^{\gamma}_n(\theta)$ on the hyperparameters $\theta$ is not known a priori, we optimize $\theta$ using an estimation of the gradient of $\mathcal{C}^{\gamma}_n$ based on \gls{DAS}. In practice, the two levels described in \cref{eq:BOA1,eq:BOA2} are iterated one after the other as depicted in \cref{fig:BOA}. To ensure feasibility, \gls{ISF} is applied on all solution candidates $s_n \in \mathcal{S}_{n}$ at the end of each iteration (which are first transformed back into the \gls{QUBO} domain).

The current implementation of \acrshort{HBS} does not utilize \gls{ISI}, but it could be included in future implementations. Similarly, the Lagrangian multipliers $\lambda$ have been considerd as constants, but could in principle also be included within the list of hyperparameters tuned automatically by \gls{DAS}. Penalty terms associated with box constraints are currently not included in the Ising Hamiltonian utilized by \gls{HBS}.  While we tailor \gls{HBS} towards its best performance for the present use case, it can, unlike \gls{IQTS}, be easily be adapted to be applicable to any \gls{QUBO}. An exemplary demonstration of the convergence behavior of \gls{HBS} is provided in \cref{supp:hbs}.

\section{Experiments} \label{sec:experiments}

In this section, we describe our numerical experiments. We use the IonQ device \emph{Aria-1}~\cite{ionq} to execute quantum circuits, which is based on trapped ion technology and operates on \num{25} qubits with all-to-all connectivity. We access the system via AWS Braket. We also conduct quantum simulations, for which we use the (noiseless) AWS Braket universal state vector simulator \emph{SV1}. In total, we perform four types of experiments:
\begin{experimentlist}
    \item \label{exp:1} \Gls{IQTS} for different scalarization weights using \gls{QAOA} on \emph{Aria-1}. 
    \item \label{exp:6} \Gls{IQTS} for the same scalarization weights as in \cref{exp:1} using \gls{SA}. 
    \item \label{exp:2} \Gls{IQTS} for different scalarization weights using \gls{QAOA} on \emph{SV1}. 
    \item \label{exp:4} \Gls{HBS} for different scalarization weights using \gls{CACm}. 
    \item \label{exp:5} \Gls{HBS} for different scalarization weights using \gls{IBP}, \gls{CACm}, and \gls{QAOA} on \emph{SV1}.
    \item \label{exp:3} \Gls{IQTS} for different primary source shares using \gls{QAOA} on \emph{SV1}.
\end{experimentlist}
A summary of all experiments is listed in \Cref{tab:experiments}. Detailed settings are provided in \cref{supp:experiments}.

\begin{table}[htb!]
\centering
\caption{Experiment summary: solver, number of instances and resulting objective hypervolume over all instances. The hypervolume is calculated with respect to the reference point $(\num{3}, \num{5}, \num{4.5}, \num{5.5})$ for the four weighted \glspl{KPI} in \cref{eqn:qubo}.}\label{tab:experiments}
\begin{tabular}{cccc}
\toprule
Experiment & Solver & Instances & Hypervolume \\
\midrule
\Cref{exp:1} & \Gls{IQTS} & \num{8} & \num{55.20} \\
\Cref{exp:6} & \Gls{IQTS} & \num{8} & \num{54.44} \\
\Cref{exp:2} & \Gls{IQTS} & \num{286} & \num{73.05} \\
\Cref{exp:4} & \Gls{HBS} & \num{206} & \num{58.74} \\
\Cref{exp:5} & \Gls{HBS} & \num{156} & \num{63.90} \\
\Cref{exp:3} & \Gls{IQTS} & \num{20} & \num{53.85} \\
\bottomrule
\end{tabular}
\end{table}

We particularly choose primary source shares of $\alpha := \alpha_i := \num{0.8}$ for all $i \in I$ in \cref{exp:1,exp:6,exp:2,exp:4,exp:5}, which is a suitable value for real-world applications. In \cref{exp:3}, on the other hand we randomly sample values for $\alpha \in [.5,.8]$. Based on these specifications, we end up with $N_z := \num{494}$ ancilla variables, leading to a total of $N_x := N_y + N_z = \num{2416}$ variables in \cref{exp:1,exp:6,exp:2,exp:4,exp:5} (while the number of variables for \cref{exp:3} depends on the sampled value of $\alpha$, which may affect the number of ancilla variables). This number can be compared to the number of variables a naive approach to the assignment problem to get a feeling for the complexity reduction we achieved with our modeling approach. If we take into account all \num{48} parts, \num{43} manufacturing sites, and \num{29} suppliers of the problem instance, we arrive at around $(43 \times 29)^{48 \times 2} \approx 10^{300}$ possible assignments, when taking double sourcing into account. In addition, the instance includes \num{28} warehouses and \num{17073} transportation methods between sites and warehouses, resulting in a large number of possible routes to which the products must also be properly assigned. Our model captures the entire problem with only $\num{2416}$ binary variables.

In the following, we first we present the results from \cref{exp:1}, which we subsequently compare with the results from \cref{exp:6}. Next, we show the resulting solution space exploration from \cref{exp:2,exp:4,exp:5}. Finally, we provide the results from \cref{exp:3}.

\subsection{Experiment E1 on Aria-1}

\begin{figure}[!tb]
	\centering
    \includegraphics[]{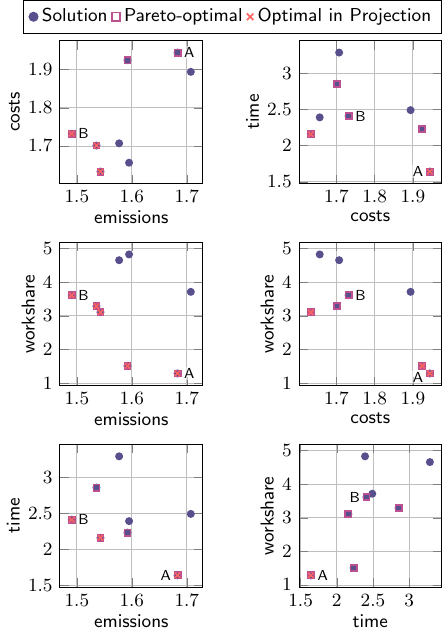}
	\caption{Solutions from \cref{exp:1}, shown as pairwise projections of the four \glspl{KPI}. Pareto-optimal solutions are highlighted. We mark two solutions, A and B, which are visualized in \cref{fig:examplesolutions}. Each solution represents a supply chain configuration.}
	\label{fig:pareto-e1}
\end{figure}

The solutions from \cref{exp:1} are shown in \cref{fig:pareto-e1} as pairwise projections within the four-dimensional \gls{KPI} space. Each solution is feasible (\ie, it fulfills the constraints) and has been obtained by choosing a different set of scalarization weights and running \gls{IQTS}. The reported \glspl{KPI} for each solution are to be understood as the non-scalarized expressions $C_n$ from \cref{eqn:C1234b}, which we refer to as costs ($n=1$), emissions ($n=2$), time ($n=3$), and workshare ($n=4$). For each two-dimensional projection, we highlight both the overall Pareto-optimal solutions as well as the solutions that are Pareto-optimal with respect to the projected plane (\ie, ignoring the other dimensions).

We also mark two Pareto-optimal example solutions, A and B, which are visualized in \cref{fig:examplesolutions}. Each visualization shows  Below the map, two plots show the workshare for the individual sites and suppliers, respectively, together with the required upper and lower bounds (triangles), which are fulfilled. For the suppliers, we also show the target workshares (crosses), which have to be reached as close as possible.

In \cref{fig:examplesolutions}a and \cref{fig:examplesolutions}b, we show an abstract map with production sites (dots) and warehouses (squares) in the style of \cref{fig:problemSketch}b. The connecting lines show the streams of goods: the line thickness indicates the cumulative costs value of the transported parts (a bigger line means more expensive), whereas the redness indicates the contribution of this connection to the overall carbon dioxide emission (the redder, the more emissions). In \cref{fig:examplesolutions}c and \cref{fig:examplesolutions}d, we show the workshare fulfillments for the sites and suppliers. The sites and suppliers are are both shown on the same horizontal axis, where each supplier is marked (triangle). Also shown are the required upper and lower limits, which are all respected by the effective workshares from the solutions (the horizontal order of sites and suppliers is with ascending feasible window size). For the suppliers, we also show the target workshares (crosses), which are to be reached as close as possible.

The visualizations demonstrate the complexity behind each individual solution and could in a next step be used to select between different competing supply chain configurations that focus on different \glspl{KPI}. For example, the solution A fulfills the supplier target workshares much better than solution solution B at a cost of a much higher carbon dioxide emission.

\ifthenelse{\boolean{ejorflag}}{%
\begin{figure*}[!t]
    \centering
    \begin{subfigure}[t]{0.3\textwidth}
        \centering
        \includegraphics[scale=1]{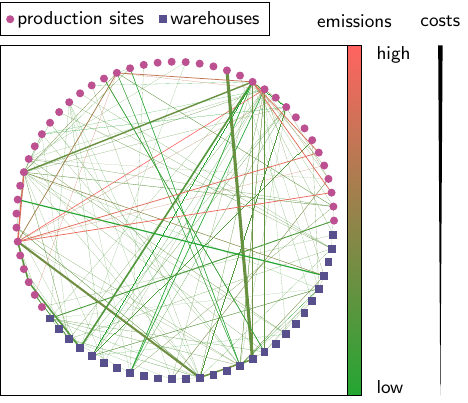}
        \caption{Solution A: Routes with emissions and costs}
        \label{fig:examplesolutions:a}
    \end{subfigure}
    \hspace{4.5cm}
    \begin{subfigure}[t]{0.3\textwidth}
        \centering
        \includegraphics[scale=1]{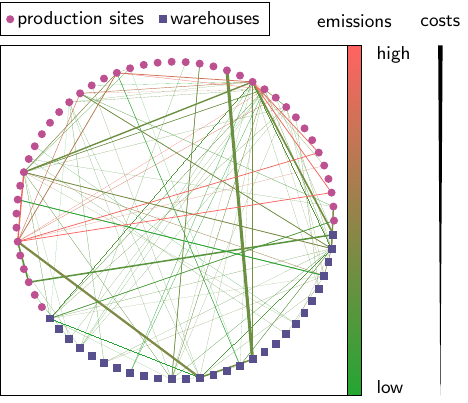}
        \caption{Solution B: Routes with emissions and costs}
        \label{fig:examplesolutions:b}
    \end{subfigure}
    \\[.2cm]
    \begin{subfigure}[t]{0.3\textwidth}
        \centering
        \includegraphics[scale=1]{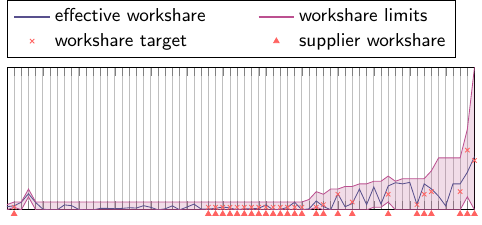}
        \caption{Solution A: Workshare fulfillments}
        \label{fig:examplesolutions:c}
    \end{subfigure}
    \hspace{4.5cm}
    \begin{subfigure}[t]{0.3\textwidth}
        \centering
        \includegraphics[scale=1]{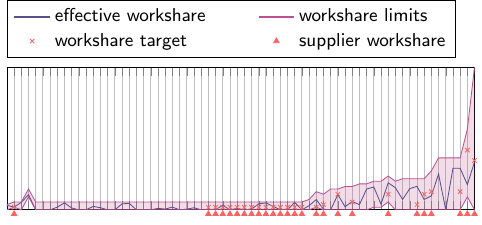}
        \caption{Solution B: Workshare fulfillments}
        \label{fig:examplesolutions:d}
    \end{subfigure}
    \caption{Visualization of two Pareto-optimal solutions from \cref{exp:1}, which we call solution A and solution B. In (a) and (b), we show routes between production sites and warehouses, the thickness and color indicate cumulative costs and emissions, respectively. In (c) and (d), we show workshare fulfillments for sites and suppliers (both on the same horizontal axis, suppliers are marked). The \glspl{KPI} of the two solutions are shown in \cref{fig:pareto-e1}.}
    \label{fig:examplesolutions}
\end{figure*}
}{%
\begin{figure*}[!t]
    \centering
    \begin{subfigure}[t]{0.49\textwidth}
        \centering
        \includegraphics[scale=1]{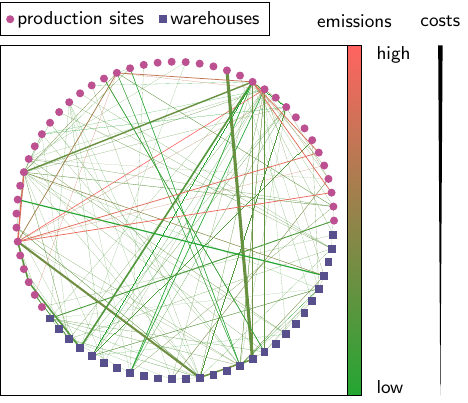}
        \caption{Solution A: Routes with emissions and costs}
        \label{fig:examplesolutions:a}
    \end{subfigure}
    ~
    \begin{subfigure}[t]{0.49\textwidth}
        \centering
        \includegraphics[scale=1]{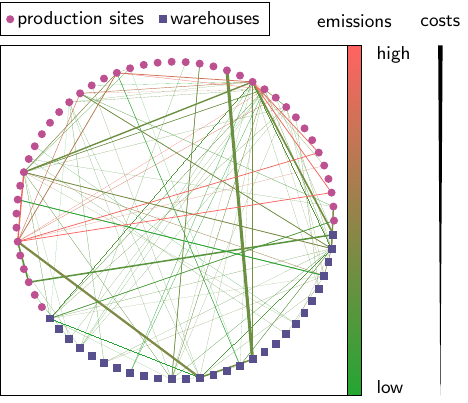}
        \caption{Solution B: Routes with emissions and costs}
        \label{fig:examplesolutions:b}
    \end{subfigure}
    \\[.2cm]
    \begin{subfigure}[t]{0.49\textwidth}
        \centering
        \includegraphics[scale=1]{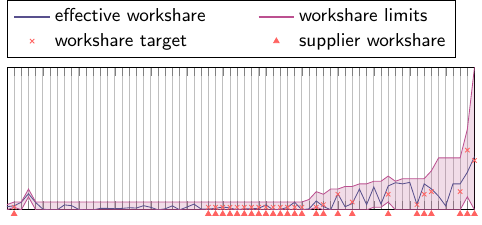}
        \caption{Solution A: Workshare fulfillments}
        \label{fig:examplesolutions:c}
    \end{subfigure}
    ~
    \begin{subfigure}[t]{0.49\textwidth}
        \centering
        \includegraphics[scale=1]{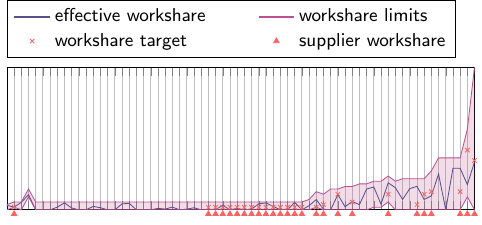}
        \caption{Solution B: Workshare fulfillments}
        \label{fig:examplesolutions:d}
    \end{subfigure}
    \caption{Visualization of two Pareto-optimal solutions from \cref{exp:1}, which we call solution A and solution B. In (a) and (b), we show routes between production sites and warehouses, the thickness and color indicate cumulative costs and emissions, respectively. In (c) and (d), we show workshare fulfillments for sites and suppliers (both on the same horizontal axis, suppliers are marked). The \glspl{KPI} of the two solutions are shown in \cref{fig:pareto-e1}.}
    \label{fig:examplesolutions}
\end{figure*}
}%

\begin{figure}[!htb]
	\centering
    \includegraphics[]{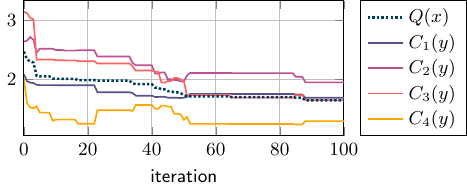}
	\caption{Convergence of \gls{IQTS} for one instance within \cref{exp:1} using \gls{QAOA} on \emph{Aria-1}.}
	\label{fig:ara1plot}
\end{figure}

In \cref{fig:ara1plot}, we present the iterative convergence of \gls{IQTS} for one instance from \cref{exp:1} (which leads to solution A). Here, $Q(x)$ stands for the objective from \cref{eqn:qubo}, which is chosen as the weighted mean of the four \glspl{KPI} from \cref{eqn:C1234}: carbon dioxide emissions ($C_1(y)$), costs ($C_2(y)$), time ($C_3(y)$), and supplier workshare target fulfillment ($C_4(y)$). It turns out that after \num{50} iterations, convergence is almost reached, which is why we use this as the maximum iteration count for \cref{exp:3}.

\subsection{Experiment E2 with SA}

The quantum device we use is capable of solving sub-problems involving up to \num{25} optimization variables using the proposed \gls{QAOA} approach. \Glspl{QUBO} of this size can also be solved relatively easy with a brute-force methods. Therefore, no concrete advantage is expected from applying a quantum optimization strategy to such small instances. However, solving sub-problems to optimality does not necessarily yield the best overall performance. In fact, sub-optimal solutions to individual sub-problems may contribute to a better overall outcome. Consequently, the inherent imperfections of the quantum optimization process may prove beneficial to the performance of \gls{IQTS}.

\begin{figure}[!tb]
	\centering
    \includegraphics[]{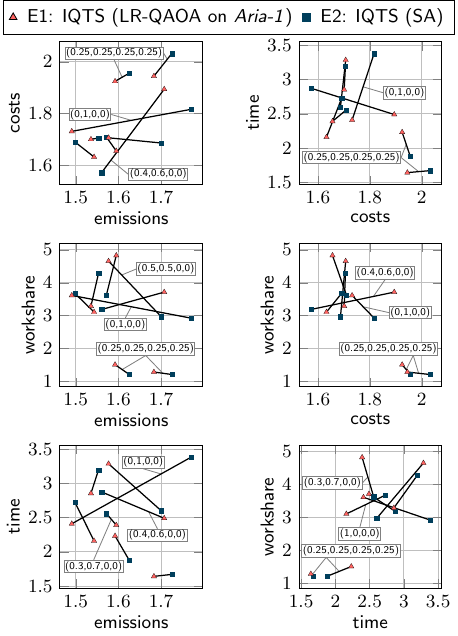}
	\caption{Comparison of the solutions from \cref{exp:1,exp:6}. Solution pairs from instances with the same scalarization weights $w$ are connected with a line. Hence, each line shows the \gls{KPI} shifts from using \gls{IQTS} with \gls{QAOA} against \gls{IQTS} with \gls{SA}. For some solution pairs, the corresponding choice of $w$ is shown.}
	\label{fig:pareto-iqts}
\end{figure}

In \cref{exp:6}, we replace \gls{QAOA} as a sub-solver with \gls{SA} to test if any performance difference can be observed. We use the same scalarization weights for the \num{8} instances as in \cref{exp:1}. A direct comparison of the results is shown in \ref{fig:pareto-iqts}, where we particularly highlight the \gls{KPI} shifts of solutions with the same scalarization weights.

There occur significant \gls{KPI} shifts, which is no surprise as the feasible solution space is very complex and many compromises may lead to similar overall objectives. We observe that the overall solution space exploration is very similar for \cref{exp:1,exp:2}. In \cref{tab:experiments}, we evaluate the hypervolume of the solutions with respect to a reference point (Nadir point). The hypervolume serves as an indicator of the Pareto frontier expansion, with higher values reflecting a better performance. We find a hypervolume of \num{55.20} for \cref{exp:1} and \num{54.44} for \cref{exp:6}, indicating a comparable Pareto frontier expansion. This quantitative result aligns with the qualitative observations presented in \cref{fig:pareto-iqts}.

In summary, no significant difference between using \gls{QAOA} or \gls{SA} in \gls{IQTS} can be observed, as expected. However, it should be emphasized that we only consider \num{8} instances per experiment here. A study with more data might be necessary to statistically validate these observations and draw more robust conclusions.

\subsection{Feasible Solution Space from E3, E4, and E5}

\begin{figure}[!tb]
	\centering
    \includegraphics[]{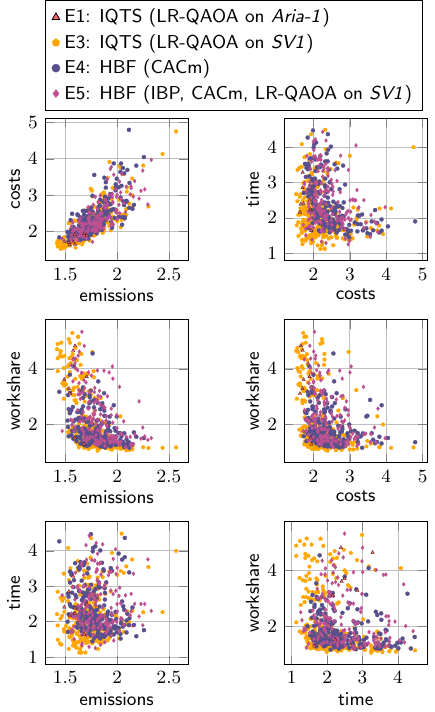}
	\caption{Explored solution space for different solvers. Each solution represents a supply chain configuration. Data from \cref{exp:1,exp:2,exp:4,exp:5}.}
	\label{fig:pareto}
\end{figure}

To visualize the solution space exploration for \cref{exp:2,exp:4,exp:5}, we show the combined results as projections on the \glspl{KPI} in \cref{fig:pareto} in analogy to \cref{fig:pareto-e1}. We also include the results from \cref{exp:1} as a reference. All presented solutions are feasible. For the sake of clarity, we do not explicitly indicate the Parto optimality of individual solutions here. These results are provided in \cref{supp:exp}.

We observe that all of our approaches lead to comparable solution space exploration with only minor difference. From \cref{tab:experiments} we find that the hypervolume of \cref{exp:2} is largest with \num{73.05} followed by \cref{exp:5} with \num{63.90} and \cref{exp:4} with \num{58.74}. However, with \num{286} instances, \cref{exp:2} also includes more solutions than \cref{exp:4} with \num{206} instances and \cref{exp:5} with only \num{156} instances. Overall, we conclude that the performance of all considered approaches is competitive. As in our previous observations, using \emph{Aria-1} in \cref{exp:1} does not yield a significant computational advantage over the classical algorithms. 

It is evident that most objectives are in conflict with one another, such as costs-time, emissions-workshare, costs-workshare, emissions-time, and time-workshare. In contrast, emissions and costs are hardly conflicting. This observation suggests that optimized supply chains can be environmentally friendly and cost-effective at the same time.

\subsection{Source Share Variation from E6}

In contrast to the previous experiments, the goal of \cref{exp:3} is to explore the effect of different values of the primary source share $\alpha$ (for constant, balanced scalarization weights $w$).

\begin{figure}[!htb]
	\centering
    \includegraphics[]{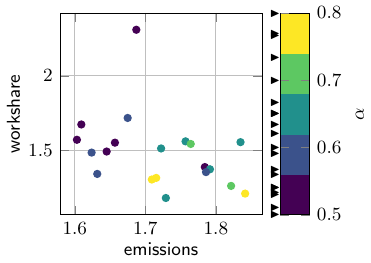}
	\caption{Solutions for different primary source shares $\alpha$. Data from \cref{exp:3}.}
	\label{fig:alphas}
\end{figure}

We present the resulting solutions in \cref{fig:alphas} for one pair of \glspl{KPI}, the other projections are shown in \cref{supp:exp}. The colors of the dots indicate the corresponding choice of $\alpha$, grouped into five distinct clusters, as defined in the colorbar. The marks on the colorbar indicate the sampled values of $\alpha$, each of which leads to a problem instance to be solved.

As expected, different choices of $\alpha$ lead to a distribution of solutions. While there is no clear distinction, higher values of $\alpha$ have a slight tendency to lead to solutions with higher emissions and lower workshare (and vice-versa). This aligns with the intuitive understanding that higher values of $\alpha$ together with the workshare constraints reduce the flexibility of choosing emission-efficient routes.

\section{Conclusions} \label{sec:conclusions}

We address a multi-objective logistics optimization problem inspired by a real-world supply chain scenario. To model the assignment problem in a form suitable for quantum optimization, we first formulate a \gls{QUBO}. By preprocessing the instance data through pathfinding and feasibility exploration, we substantially reduce the solution space that must be searched. Building on this foundation, we propose two hybrid quantum-classical solvers to tackle the optimization problem effectively: The knowledge-based \gls{IQTS} and the high-performance approach \gls{HBS}. Both approaches leverage special-purpose systems to solve sub-problems of scalable size, enabling the use of Ising machines, \gls{NISQ} devices, and potentially more efficient quantum or non-classical platforms in the future.

We demonstrate experimentally that the two heuristic solvers are able to find solutions of comparable quality. Both have their own strengths and weaknesses. \Gls{IQTS} is highly efficient in finding feasible solutions and its focus on the exploration of the solution space close to the feasible domain typically leads to a very fast convergence behavior. Its weakness, on the other hand, is that it cannot be directly parallelized in its current design. The strength of \gls{HBS} is its flexibility in combining multiple sub-solvers within a bilevel framework, which ensures scalability and generalization. A weakness of the method is that its adaptability to a general problem setting comes at the cost of slower convergence to an optimal solution. In the course of our experiments, however, the two solvers demonstrate similar solution quality, suggesting that these weaknesses were minor in practice.

Our experimental results do not indicate any performance advantage from using the \emph{Aria-1} quantum hardware within the proposed solver frameworks. We attribute this primarily to the small scale of the device of only 25 qubits, which cannot be expected to be competitive to classical computations~\cite{dalzell2020,xu2023}. Nevertheless, the design of our solver frameworks allows for straightforward adaptation to more powerful quantum hardware in the future, where quantum optimization may lead to meaningful performance gains.


There are numerous opportunities to extend and generalize our modeling approach from multiple perspectives. One important aspect would be the packing problem of placing multiple parts within the limited container space of transport vehicles based on the geometries. Additional elements could also be taken into account, such as production times, precise transportation scheduling, manufacturing capacities, warehouse capacities, and additional regional restrictions related to transportation and regulatory compliance. Some of these elements, such as production times, could be easily integrated into our existing formalism. Others, such as the packing problem, would lead to entirely new sub-problems that could be integrated as a future extensions. Lastly, the primary source share, which we treated an instance-dependent example, could also be taken into account as an additional optimization parameter. Possibly, in combination with another \gls{KPI}, which measures the resilience of the supply chain.

In the current work, we have not combined all elements of \gls{IQTS} and \gls{HBS}, as integrating calls to the quantum hardware with parallel CPU processing is non-trivial in the current setup. However, this integration could be achieved with improved interfacing. \Gls{HBS} can be extended to include a broader set of special-purpose algorithms, either classical, quantum-inspired or quantum. Incorporating state-of-the-art algorithms into the solver mix ensures that the method achieves optimal performance at present. As quantum hardware scales to include more qubits, it becomes possible to allocate an increasing portion of the solver mix to quantum algorithms. In the meantime, the hybrid approach facilitates the efficient utilization of these schemes. Moreover, \gls{HBS} can be improved by incorporating more sophisticated methods for merging candidate solutions from solvers (such as mutation rules) and by exploring alternative parameter tuning techniques, like Bayesian optimization, where applicable. For classical hardware, this study has primarily focused on CPU-based parallel computing, as the sparsity of the problem's connectivity made CPUs a suitable choice. However, algorithms such as \gls{CACm} and \gls{IBP} could benefit from GPU acceleration in cases of denser connectivity.

Another potential non-classical platform for Ising machines, specifically of the \gls{CACm} type, is the \gls{CIM}. The \gls{CIM} is a nonlinear optical system that leverages optical parametric amplification which can be implemented on the lithium niobate platform~\cite{honjo2021100,Marandi2014,mcmahon2016fully,inagaki2016large}. A next-generation \gls{CIM} should be implemented as a \gls{TFLN} photonic integrated circuit. Several key components toward this goal have successfully developed such as optical parametric amplifiers~\cite{jankowski2022quasi}, optical frequency converters~\cite{jankowski2023supercontinuum} and optical parametric oscillators~\cite{mckenna2022ultra}.

\section{Acknowledgments}
The authors acknowledge Airbus and BMW Group for providing the problem statements under the Airbus–BMW Group Quantum Computing Challenge 2024 and for supporting the validation of the results. We also extend our appreciation to the teams at Airbus, Amazon, and IonQ for their invaluable guidance in clarifying the use case and for providing access to IonQ’s computing resources via AWS. Their readiness to address our questions throughout the development process was especially appreciated.
We would like to give special thanks (in alphabetical order) to Andreas Mitschke (Airbus), Andreas Zindel (Airbus), Denny Dahl (IonQ), Heidi Nelson-Quillin (IonQ), Howard Lock (IonQ), Martin Schuetz (Amazon), Michael Brett (Amazon), Richard Ashworth (Airbus), and Sebastian Stern (Amazon). Furthermore, we express our sincere gratitude to Marwan Channab, Raymond Harding, and Matteo Conterno for their active participation in discussions and the exchange of ideas.

\printbibliography

\begin{appendices}
\section{Preprocessing} \label{app:preprocessing}

We consider here \cref{eqn:model1} with scalarized objectives $\sum_{n=1}^4 w_n C_n(y)$. As described in \cref{sec:qubo model}, we perform a two-step preprocessing step to remove sub-optimal terminals from the formulation and reduce the overall number of variables. These two prepreocessing steps are explained in the following.

\subsection{Pathfinding}

So far, we have not specified how to determine the set of all routes from the set of transportation methods, which is necessary to define the terminals for each site. It is clear that an explicit enumeration of all routes is computationally expensive due to the many combinatorial possibilities. However, it turns out that with our scalarization approach, no explicit representation is necessary because it allows us to reduce the number of routes for each part and each pair of sites to only one, as will be explained in the following.

By definition, the first three \glspl{KPI}, \cref{eqn:C123b}, consist of a sum of non-negative contributions from the transport costs, \cref{eqn:cover}. Therefore, when omitting the contribution from $C_4(y)$, the scalarized objective reads
\begin{align} \label{eqn:C123w}
    & \sum_{n \in \{1,2,3\}} w_n C_n(y) \\
    = & \sum_{\substack{(i,j) \in \phi \\(k,t,u) \in F_i,a \in A \\ (l,s,v) \in F_j, b \in A \\ R^i_{k \Rightarrow l} \neq 0}} \!\!\!\!\!\!  \alpha_i^a \chi^{i}_{kl} (\overline{w},t,s) \!\!\!\!\! 
    & y_{i\rightarrow(k,t.u)}^{a} y_{j\rightarrow(l,s,v)}^{b} \nonumber,
\end{align}
where
\begin{align} \label{eqn:chi}
    \chi^{i}_{kl}(\overline{w},t,s) := \sum_{\substack{n\in\{1,2,3\}}} \frac{w_n}{d_n} c_{(k,t)\Rightarrow(l,s)}^{in}
\end{align}
represents the contribution of transporting part $i \in I$ from production site $k \in K$ and terminal $t \in T^i_k$ to production site $l \in K$ and terminal $s \in T^j_l$. The end terminal $s$ has in fact no influence on \cref{eqn:chi}, so that we can also write $\chi^{i}_{kl}(\overline{w},t) := \chi^{i}_{kl}(\overline{w},t,s)$. Moreover, \cref{eqn:chi} and hence \cref{eqn:C123w} is independent of the chosen suppliers. We use $\overline{w} := (w_1,w_2,w_3)$ to denote the first three scalarization weights.

The start terminal $t$ represents the chosen route $r = \tau^i_{kl}(t) \in R^i_{k \Rightarrow l}$ (as a sequence of transportation options) for part $i$. However, since every part is to be transported independently, there is in fact a universally optimal route for each part $i$ between each pair of production sites $(k,l)$, depending on the chosen weights $\overline{w}$, which is formally given by the optimal terminal 
\begin{align} \label{eqn:bestt}
    \hat{t}^{i}_{kl}(\overline{w}) := \argmin_{t \in T^i_{k}} \chi^{i}_{kl}(\overline{w},t).
\end{align}
This optimal terminal can in principle be found by enumerating all possible terminals $T^i_k$. Alternatively, finding the best terminal is equivalent to finding the most cost-efficient path between $k$ and $l$ on the graph $G_i$, where each edge $m \in M_i$ is assigned the costs
\begin{align} \label{eqn:Gcost}
    \xi^i(\overline{w},m) := \sum_{\substack{n\in\{1,2,3\}}} \frac{w_n}{d_n} c^n(m) \gamma^{in}(m)
\end{align}
since
\begin{align}
    \chi^{i}_{kl}(\overline{w},t) = \sum_{m \in \tau^i_{kl}(t)} \xi^i(\overline{w},m),
\end{align}
where we recall \cref{eqn:cover,eqn:chi}. $R^i_{k \Rightarrow l} \neq \emptyset$, the most cost-efficient path on $G_i$ corresponds to a route
\begin{align} \label{eqn:bestr}
    \hat{r}^i_{kl}(\overline{w}) := \argmin_{r \in R^i_{k \Rightarrow l}} \sum_{m \in r} \xi^i(\overline{w},m)
\end{align}
with
\begin{align}
    \hat{t}^{i}_{kl}(\overline{w}) = \overline{\tau}^i_{kl}(\hat{r}^i_{kl}(\overline{w})).
\end{align}
For $R^i_{k \Rightarrow l} = \emptyset$, the optimal terminal is simply the only available terminal, \cref{eqn:numterminals}.

We consider the determination of optimal terminals, as described above, as a preprocessing step, which we perform for each choice of weights $\overline{w}$. To this end, \cref{eqn:bestr} is solved with Dijkstra's algorithm. This pathfinding step eliminates the search for the best terminal (or, equivalently, best route) from the optimization problem and allows us to switch to a new set of of binary variables $\overline{y}$, each identifiable by four indices. Specifically, the binary variable $\overline{y}_{i\rightarrow(k,u)}^{a} \in \{0,1\}$ represents the assignment of part $i \in I$ to the best terminal $\hat{t}^{i}_{kl}(\overline{w}) \in T^i_k$ of site $k$ (presuming that parent part $j in I$ is assigned to production site $l \in K$) and supplier $u$ with $(k,u) \in f_i$, where $a \in A$ indicates either the primary production ($a=1$) or the secondary production ($a=2$). Furthermore, we write
\begin{align}
    c^{in}_{k \Rightarrow l} (w_1,w_2,w_3) := c^{in}_{(k,\hat{t}^{i}_{kl}(\overline{w}))\Rightarrow (l,s)}
\end{align}
for the cumulative contributions of the transportation methods, where we recall the independence from the end terminal $s$.

\subsection{Feasibility Reduction}

In addition to the pathfinding, we apply a second preprocessing step in which we refine the feasible solution space. This step is independent of the chosen weights and only has to be performed once per problem instance. Specifically, for each part $i \in I$, we consider only those elements from the predefined set of feasible site-supplier combinations $f_i$ that contain sites for which every child of $i$ can be transported to $i$ (as its parent) and, at the same time, $i$ can also be transported to its parent. Given these conditions, the resulting reduced set is given by the self-referencing definition
\begin{align}
    g_i := & \{ (k, u)  \,|\, (k,u) \in f_i \land [ \forall (i',j) \in \phi \\
    & \land i'=i \exists (l, v) \in g_j : R^i_{k \Rightarrow l} \neq \emptyset ] \nonumber \\
    & \land [\forall (j,i') \in \phi \land i'=i \exists (l, v) \in g_j \nonumber \\
    & : R^j_{l \Rightarrow k} \neq \emptyset] \} \subseteq f_i \nonumber
\end{align}
and can be straightforwardly determined with a recursive search of $f_i$. We find that $\sum_{i \in I}|g_i| = \num{1922}$. Since all options from $f_i \setminus g_i$ violate \cref{eqn:model1:cI}, we only use the choices provided by $g_i$ in the following, leading to a total of $N_y := \num{1922}$ binary variables. Furthermore, $g_i$ also allows us to fix assignments with a single feasible option, \ie, $\overline{y}_{i\rightarrow(k,u)}^{a} := 1$ for all $i \in I$ with $|g_i|=1$. However, this approach does not lead to a reduction of variables for the problem instance we consider.

\subsection{Reduced Model}

As a consequence of the two preprocessing steps, the reduced model is given by:

\begin{condstrip}
\begin{subequations}
\begin{align}
\min_{\overline{y}} & \sum_{n=1}^4 \overline{C}_n(w,\overline{y}) \label{eqn:model2} \tag{\theparentequation} \\
    \text{s.t. } & \sum_{a,b \in A, u', v' \in U} \mu_i^{ab} \!\!\! \delta_{u'u} \delta_{v'v} \overline{y}_{i\rightarrow(k,u)}^{a} \overline{y}_{j\rightarrow(l,v)}^{b} = 0 \,\forall\, (i,j) \in \phi \,\forall\, (k, u) \in g_i \,\forall\, (l, v) \in g_j \land R^i_{k \Rightarrow l} = \emptyset \label{eqn:model2:cI} \\
     & \sum_{(k,u) \in g_i} \overline{y}_{i\rightarrow(k,u)}^{a} = 1 \,\forall\, i \in I \,\forall\, a \in A \label{eqn:model2:cII} \\
     & \Big( \sum_{(k',u) \in g_i} \delta_{k'k} \overline{y}_{i\rightarrow(k,u)}^{a=1} \Big) \Big( \sum_{(k',u) \in g_i} \delta_{k'k} \overline{y}_{i\rightarrow(k,u)}^{a=2} \Big) = 0 \,\forall\, i \in I \,\forall\, k \in K \land |K_i|\geq2 \label{eqn:model2:cIII} \\
     & \Big( \sum_{(k,u) \in g_i} V_{vk} \overline{y}_{i\rightarrow(k,u)}^{a=1} \Big) \Big( \sum_{(k,u) \in g_i} V_{vk} \overline{y}_{i\rightarrow(k,u)}^{a=2} \Big) = 0 \,\forall\, i \in I \,\forall\, v \in V \land |V_i|\geq2 \label{eqn:model2:cIV} \\
     & \sum_{i \in I, (k',u) \in g_i, a \in A} \delta_{k'k} v_i \alpha_i^a \overline{y}_{i\rightarrow(k,u)}^{a} \in [K^{\text{min}}_{k},K^{\text{max}}_{k}] \,\forall\, k \in K \label{eqn:model2:cV} \\
     & \sum_{i \in I, (k,u') \in g_i, a \in A} \delta_{u'u} v_i \alpha_i^a \overline{y}_{i\rightarrow(k,u)}^{a} \in [U^{\text{min}}_{u},U^{\text{max}}_{u}] \,\forall\, u \in U \label{eqn:model2:cVI}
\end{align}
\end{subequations}
\end{condstrip}

Here, we use the simplified \glspl{KPI}
\begin{subequations} \label{eqn:C1234}
\begin{align} \label{eqn:C123}
    & \overline{C}_n(w,\overline{y}) := \overline{C}_n(w_1,w_2,w_3,\overline{y}) \\
    & := \frac{w_n}{d_n} \sum_{\substack{(i,j) \in \phi \\(k,u) \in g_i,a \in A \\ (l,v) \in g_j, b \in A \\ R^i_{k \Rightarrow l} \neq 0}} \alpha_i^a c^{in}_{k \Rightarrow l} (w_1,w_2,w_3) \overline{y}_{i\rightarrow(k,u)}^{a} \overline{y}_{j\rightarrow(l,v)}^{b} \nonumber
\end{align}
for $n \in \{1,2,3\}$ and
\begin{align} \label{eqn:C4}
    & \overline{C}_4(w,\overline{y}) := \overline{C}_4(w_4,\overline{y}) \\
    & := \frac{w_4}{d_4} \sum_{u \in U} \Big[ \!\!\!\! \sum_{\substack{i \in I\\(k,v) \in g_i\\a \in A}} \!\!\! \delta_{uv} v_i \alpha_i^a \overline{y}_{i\rightarrow(k,u)}^{a} - U^{\text{target}}_{u} \Big]^2 \nonumber,
\end{align}   
\end{subequations}
respectively, in analogy to \cref{eqn:C123,eqn:C4}.

\section{Box Constraints} \label{app:box constraints}

As described in \cref{sec:qubo model}, the box constraints from \cref{eqn:model2}, \cref{eqn:model2:cV,eqn:model2:cVI}, must first be converted into equality constraints before they can be transformed into penalty terms. We describe the procedure in the following~\cite{babbush2013}.

\Cref{eqn:model2:cV,eqn:model2:cVI} can also be written in terms of inequalities $f(\overline{y}) \geq 0$ using functions of the form
\begin{align}
    f(\overline{y}) := \sum_{i \in I,(k,u) \in g_i,a \in A} b^i_{(k,u)} v_i \alpha_i^a \overline{y}^a_{i \rightarrow (k,u)} + \nu
\end{align}
with coefficients $b^i_{(k,u)} \in \{-1,0,1\}$ for $(k,u) \in g_i$ and an integer $\nu \in \mathbb{Z}$.

As an initial step, and for reasons that will become clear further below, we approximate the real-valued coefficients (the relative values and the global primary source shares) by rational numbers. To that end, we define a part-dependent numerators $P_i,\overline{P}_i \in \mathbb{Z}$ for each part $i \in I$ and global denominators $R,\overline{R} \in \mathbb{Z}_{\neq0}$, which allow us to express the relative value as
\begin{align} \label{eqn:vapprox}
    v_i = \frac{P_i}{R} + \epsilon_i
\end{align}
and the primary source share as 
\begin{align} \label{eqn:aapprox}
    \alpha_i = \frac{\overline{P}_i}{\overline{R}} + \overline{\epsilon}_i,
\end{align}
respectively, with error terms $\epsilon_i,\overline{\epsilon}_i \in \mathbb{R}$ for all $i \in I$. Thus, we can write
\begin{align}
    f(\overline{y}) = f_{\mathbb{Q}}(\overline{y}) + f_{\epsilon}(\overline{y})
\end{align}
with a rational-valued part
\begin{align}
    f_{\mathbb{Q}}(\overline{y}) := \frac{f_{\mathbb{Z}}(\overline{y})}{R \overline{R}} \in \mathbb{Q}
\end{align}
and a real-valued part
\begin{align}
    f_{\epsilon}(\overline{y}) := \!\!\!\! \sum_{\substack{i \in I\\(k,u) \in g_i\\a \in A}} \!\! \left[ \frac{P_i}{R}\overline{\epsilon}_i + \frac{\overline{P}_i}{\overline{R}} \epsilon_i + \epsilon_i \overline{\epsilon}_i \right] \overline{y}^a_{i \rightarrow (k,u)} \in \mathbb{R}.
\end{align}
The rational-valued part is based on the integer-valued expression
\begin{align}
    f_{\mathbb{Z}}(\overline{y}) := \!\!\!\! \sum_{\substack{i \in I\\(k,u) \in g_i\\a \in A}} \!\! b^i_{(k,u)} P_i \overline{P}_i \overline{y}^a_{i \rightarrow (k,u)} + \nu \in \mathbb{Z}.
\end{align}
In the following, we presume that the numerators and denominators are chosen in such a way that $f_{\epsilon}(\overline{y}) \approx 0$ for all $\overline{y} \in \{0,1\}^{N_y}$, hence $f(\overline{y}) \approx f_{\mathbb{Q}}(\overline{y})$.

In the next step, we define the expression 
\begin{align}
    g_{\mathbb{Z}}(\overline{y},z) := f_{\mathbb{Z}}(\overline{y}) - \sum_{n=1}^{n_z} 2^{n-1} z_n \in \mathbb{Z}.
\end{align}
Here, $z_n \in \{0,1\}$ for $n\in \{1,\dots,n_z\}$ denote the $n_z$ ancilla variables with $n_z := \lceil \log_2 \overline{f}_{\mathbb{Z}} \rceil$ with a constant $\overline{f}_{\mathbb{Z}} \in \mathbb{Z}$, which can be freely chosen under the requirement $\overline{f}_{\mathbb{Z}} \geq 1+\max_{\overline{y}} f_{\mathbb{Z}}(\overline{y})$. As a consequence, $\forall \overline{y} \in \{0,1\}^{N_y} \,\exists\, z : g_{\mathbb{Z}}(\overline{y},z) = 0$ (since $f_{\mathbb{Z}}(\overline{y}) \in [0,\overline{f}_{\mathbb{Z}}]$) and we can replace the inequality constraint $f_{\mathbb{Q}}(\overline{y}) \geq 0$ by the equality constraint $g_{\mathbb{Z}}(\overline{y},z)=0$. We use this approach to replace all box constraints, \cref{eqn:model2:cV,eqn:model2:cVI}, by equality constraints. 

By definition, the minimal number of ancilla variables depends on the magnitude of the global denominators $R$ and $\overline{R}$ such that smaller values lead to less variables. More ancilla variables make the problem harder to solve. Therefore, the rational approximation in \cref{eqn:vapprox,eqn:aapprox} needs to be balanced between sufficiently small error terms and a sufficiently small number of ancilla variables. In the following, we omit the dependency on the specific choices of numerators and denominators to simplify the notation. Additionally, we aggregate all ancilla variables which emerge from the conversion of the box constraints into a set of binary variables, denoted by $z$. We then represent the combined vector of $\overline{y}$ and $z$ as $x := (\overline{y}, z)$.

The resulting model only consists of equality constraints, which can then be transformed into penalty terms. 

\section{QUBO specifications} \label{app:qubo}

\Cref{eqn:qubo:Q} from \cref{eqn:model2} makes use of the following notations. The weighted \glspl{KPI} $\overline{C}_n(w_n,\overline{y})$ are provided in \cref{eqn:C123,eqn:C4}. The penalty terms are defined by
\begin{subequations}
\begin{align}
    P_1(\overline{y}) & := \!\!\!\!\!\! \sum_{\substack{(i,j) \in \phi \\ (k, u) \in g_i, (l, v) \in g_j \\ R^i_{k \Rightarrow l} = \emptyset}} \!\!\! \Big( \! \sum_{\substack{a,b \in A\\u', v' \in U}} \mu_i^{ab}  \\
    & \hspace{2.25cm} \times \delta_{u'u} \delta_{v'v} \overline{y}_{i\rightarrow(k,u)}^{a} \overline{y}_{j\rightarrow(l,v)}^{b} \Big), \nonumber \\
    P_2(\overline{y}) & := \sum_{\substack{i \in I, a \in A}} \Big( \sum_{(k,u) \in g_i} \overline{y}_{i\rightarrow(k,u)}^{a} - 1 \Big)^2, \\
    P_3(\overline{y}) & := \sum_{\substack{i \in I, k \in K \\ |K_i|\geq2}} \Big( \sum_{(k',u) \in g_i} \delta_{k'k} \overline{y}_{i\rightarrow(k,u)}^{a=1} \Big) \\
    & \hspace{1.8cm} \times \Big( \sum_{(k',u) \in g_i} \delta_{k'k} \overline{y}_{i\rightarrow(k,u)}^{a=2} \Big), \nonumber\\
    P_4(\overline{y}) & := \sum_{\substack{i \in I, v \in V \\ |V_i|\geq2}} \Big( \sum_{(k,u) \in g_i} V_{vk} \overline{y}_{i\rightarrow(k,u)}^{a=1} \Big) \\
    & \hspace{1.8cm} \times \Big( \sum_{(k,u) \in g_i} V_{vk} \overline{y}_{i\rightarrow(k,u)}^{a=2} \Big), \nonumber\\
    P_5(x) & := P_5^{\geq}(x) + P_5^{\leq}(x), \\
    P_6(x) & := P_6^{\geq}(x) + P_6^{\leq}(x)
\end{align}
\end{subequations}
based on the abbreviations
\begin{subequations}
\begin{align}
    P_5^{\geq}(x) & := \sum_{k \in K} 2^{-n^{(5,\geq)}_k} \Big( p_5(\overline{y}) - K^{\text{min}}_{k} R \overline{R} \\
    & \hspace{2.5cm} - \sum_{n=1}^{n^{(5,\geq)}_k} 2^{n-1} z^{(5,\geq)}_{kn} \Big)^2, \nonumber\\
    P_5^{\leq}(x) & := \sum_{k \in K} 2^{-n^{(5,\leq)}_k} \Big( K^{\text{max}}_{k} R \overline{R} - p_5(\overline{y}) \\ 
    & \hspace{2.5cm} - \sum_{n=1}^{n^{(5,\leq)}_k} 2^{n-1} z^{(5,\leq)}_{kn} \Big)^2, \nonumber\\
    P_6^{\geq}(x) &  := \sum_{u \in U} 2^{-n^{(6,\geq)}_u} \Big( p_6(\overline{y}) - U^{\text{min}}_{u} R \overline{R} \\ 
    & \hspace{2.5cm} - \sum_{n=1}^{n^{(6,\geq)}_u} 2^{n-1} z^{(6,\geq)}_{un} \Big)^2, \nonumber\\
    P_6^{\leq}(x) &  := \sum_{u \in U} 2^{-n^{(6,\leq)}_u} \Big( U^{\text{max}}_{u} R \overline{R} - p_6(\overline{y}) \\ 
    & \hspace{2.5cm} - \sum_{n=1}^{n^{(6,\leq)}_u} 2^{n-1} z^{(6,\leq)}_{un} \Big)^2. \nonumber
\end{align}
\end{subequations}
The exponential factors in the sums are empirically chosen to suitably scale the penalty terms. The expressions are based on
\begin{align}
    p_5(\overline{y}) &  := \sum_{i \in I, (k',u) \in g_i, a \in A} \delta_{k'k} P_i \overline{y}_{i\rightarrow(k,u)}^{a} \overline{D}^a_i \\
    p_6(\overline{y}) & := \sum_{i \in I, (k,u') \in g_i, a \in A} \delta_{u'u} P_i \overline{y}_{i\rightarrow(k,u)}^{a} \overline{D}^a_i
\end{align}
with
\begin{align}
    \overline{D}^a_i := \begin{cases} \overline{P}_i & , \text{if } a=1 \\ \overline{R}-\overline{P}_i & , \text{if } a=2\end{cases}
\end{align}
and the ancilla counters
\begin{subequations}
\begin{align}
    n^{(5,\geq)}_k & := \Big\lceil \log_2 \Big( p_5(\overline{y}=1) - K^{\text{min}}_{k} R \overline{R} + 1\Big) \Big\rceil, \\
    n^{(5,\leq)}_k & := \Big\lceil \log_2 \Big( K^{\text{max}}_{k} R \overline{R} + 1\Big) \Big\rceil,  \\
    n^{(6,\geq)}_u & := \Big\lceil \log_2 \Big( p_6(\overline{y}=1) - U^{\text{min}}_{u} R \overline{R} + 1\Big) \Big\rceil, \\
    n^{(6,\leq)}_u & := \Big\lceil \log_2 \Big( U^{\text{max}}_{u} R \overline{R} + 1\Big) \Big\rceil.
\end{align}
\end{subequations}
The ancilla counters can be summed up to obtain the total number of ancilla variables
\begin{align}
    N_z := & \sum_{k \in K} \Big( n^{(5,\geq)}_k + n^{(5,\leq)}_k \Big) \\
    & + \sum_{u \in U} \Big( n^{(6,\geq)}_u + n^{(6,\leq)}_u \Big), \nonumber
\end{align}
where the vector of ancilla variables is to be understood as
\begin{align}
    z := \Big( & z^{(5,\geq)}_{kn} \forall n \in \{1,\dots,n^{(5,\geq)}_k\}, k \in K, \\
           & z^{(5,\leq)}_{kn} \forall n \in \{1,\dots,n^{(5,\leq)}_k\}, k \in K, \nonumber \\
           & z^{(6,\geq)}_{un} \forall n \in \{1,\dots,n^{(6,\geq)}_u\}, u \in U, \nonumber \\
           & z^{(6,\leq)}_{un} \forall n \in \{1,\dots,n^{(6,\leq)}_u\}, u \in U \Big).       
\end{align}

In summary, the following steps are necessary to arrive at \cref{eqn:qubo} from \cref{eqn:model1}:
\begin{enumeratedense}
    \item Choose the scalarization weights $w$.
    \item Find the optimal routes for each part with Dijkstra's algorithm (depends on $w$).
    \item Refine the feasible solution space (independent of $w$), eliminating all options that cannot be realized because of a missing connectivity. This leads to \cref{eqn:model2}.
    \item Choose numerators and denominators to obtain a rational approximation of the relative values and the primary source shares for all parts.
    \item Transform inequality constraints into penalties using $N_z$ ancilla variables.
    \item Choose the corresponding Lagrange multipliers $\lambda$. This leads to \cref{eqn:qubo}.
\end{enumeratedense}

\section{HBS Demonstration} \label{supp:hbs}

\begin{figure}[!htb]
    \centering
    \includegraphics[width=1.0\linewidth]{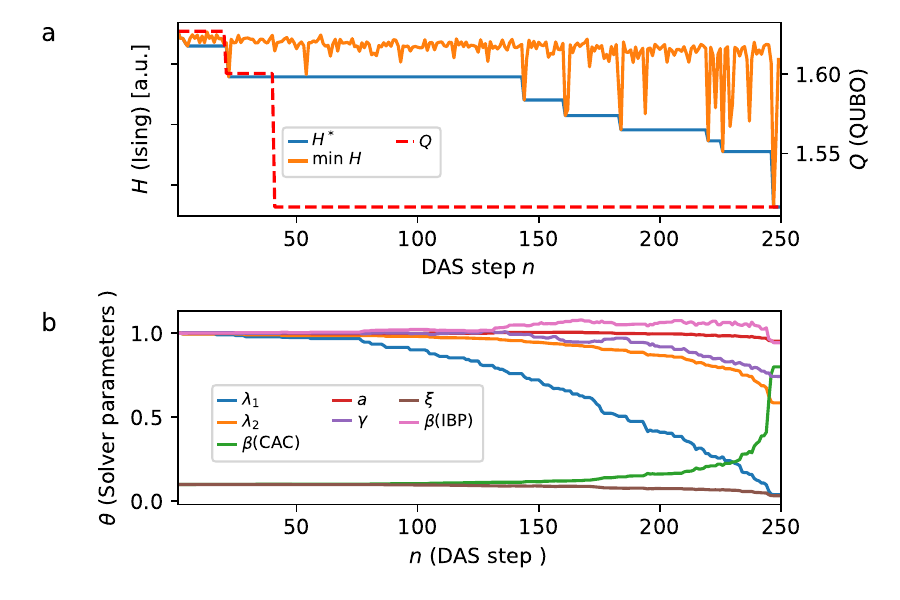}
    \caption{Simulation results of the bilevel optimization approach using \gls{CACm}, \gls{IBP}, and \gls{QAOA}. (a) Ising and \gls{QUBO} objective vs. steps of \gls{DAS}. (b) Step-dependent hyperparameters $\theta$ of the Ising solvers. $\omega_1 = 0.5$, $\omega_2 = 0.0$,$\omega_3 = 0.0$, $\omega_4 = 0.5$. See description of the hyperparameters $\{\lambda_1,\lambda_2,\gamma,\beta(\text{IBP}),\beta(\text{CAC}),\xi,a,T\}$ in \cref{sec:cacm,sec:ibp}. $H$ and $Q$ are obtained before and after the application of \gls{ISF}, respectively. The reference energy $H^*$ is defined as the minimum objective of \cref{eqn:ising} found so far.}
    \label{fig:DAS}
\end{figure}

For demonstration purposes, the time progression of the solution quality at each step $n$ of \gls{HBS} from \cref{sec:hbs} is exemplarily shown in \cref{fig:DAS}. The decrease of the Ising Hamiltonian $H$ (\cref{fig:DAS}a) is accelerated by the online tuning of hyperparameters achieved by \gls{DAS} (\cref{fig:DAS}b). At each step, the Ising solvers \gls{CACm}, \gls{IBP}, and \gls{QAOA} are applied to improve the solution quality. The solution candidates are then fixed using \gls{ISF} and their objective value is evaluated. 

\section{Experimental Setup} \label{supp:experiments}

For the experiments described in \cref{sec:experiments}, we have used the following settings (using the notation from \cref{sec:use case,sec:methods}):
\begin{itemizedense}
    \item[\ref{exp:1}.] Model: $\overline{R} := 5$, $\overline{P}_i := 4$ for all $i \in I$, $\lambda_i=2$ for $i\in\{1,6\}$, weights for each instance are chosen manually with $w_2=1-w_2$ and $w_3=w_4=0$: $w=(\num{0},\num{1},\num{0},\num{0})$, $w=(\num{0.5},\num{0.5},\num{0},\num{0})$, $w=(\num{1},\num{0},\num{0},\num{0})$, $w=(\num{0.4},\num{0.6},\num{0},\num{0})$, $w=(\num{0.6},\num{0.4},\num{0},\num{0})$, $w=(\num{0.3},\num{0.7},\num{0},\num{0})$, and twice $w=(\num{0.25},\num{0.25},\num{0.25},\num{0.25})$. \Gls{IQTS}: $m=4$, $n=15$, $\kappa=50$ ($\kappa=100$ for $w=(\num{0.25},\num{0.25},\num{0.25},\num{0.25})$). \Gls{QAOA}: $p=1$. \emph{Aria-1}: \num{256} shots per circuit. Total number of instances: \num{8}.
    \item[\ref{exp:6}.] Model, \Gls{IQTS}: as in \cref{exp:1}, but instead of \gls{QAOA}, \gls{SA} is used. \Gls{SA}: \num{1000} annealing steps. Total number of instances: \num{8}.
    \item[\ref{exp:2}.] Model: $\overline{R} := 5$, $\overline{P}_i := 4$ for all $i \in I$, $\lambda_i=2$ for $i\in\{1,6\}$, weights are chosen as a uniform grid $w_i \in \{ \num{0.1}, \num{0.2}, \num{0.3}, \num{0.4}, \num{0.5}, \num{0.6}, \num{0.7}, \num{0.8}, \num{0.9}, \num{1} \}$ for $i\in\{1,4\}$ with $\sum_{i}w_i=1$. \Gls{IQTS}: $m=4$, $n=25$, $\kappa=50$. \Gls{QAOA}: $p=1$. \emph{SV1}: \num{1024} shots per circuit. Total number of instances: \num{286}.
    \item[\ref{exp:4}.] Model: $\overline{R} := 5$, $\overline{P}_i := 4$ for all $i \in I$, $\lambda_i=3$ for $i\in\{1,4\}$, $\lambda_i=0$ for $i\in\{5,6\}$, uniformly sampled weights $w_i \in [0,1]$ for $i\in\{1,4\}$ with $\sum_{i}w_i=1$ for each instance. \Gls{DAS}: $n=250$ steps. Total number of instances: \num{206}.
    \item[\ref{exp:5}.] Model: $\overline{R} := 5$, $\overline{P}_i := 4$ for all $i \in I$, $\lambda_i=3$ for $i\in\{1,4\}$, $\lambda_i=0$ for $i\in\{5,6\}$, uniformly sampled weights $w_i \in [0,1]$ for $i\in\{1,4\}$ with $\sum_{i}w_i=1$ for each instance. \Gls{DAS}: $n=250$ steps. \Gls{QAOA}: $p=3$. \emph{SV1}: \num{32} shots per circuit. Total number of instances: \num{156}.
    \item[\ref{exp:3}.] Model: uniformly sampled denominators $\overline{R} \in \{2,3,5,10,15,20,25,30,50,100\}$ for each instance, $\overline{P}_i := \tilde{P}(\overline{R})$ for all $i \in I$ with uniformly sampled values $\tilde{P}(\overline{R}) \in \{\lceil\alpha_{\mathrm{min}}\overline{R}\rceil,\dots,\lfloor\alpha_{\mathrm{max}}\overline{R}\rfloor\}$ with $\alpha_{\mathrm{min}}:=0.5$ and $\alpha_{\mathrm{max}}:=0.8$ for each instance, $\lambda_i=2$ for $i\in\{1,6\}$, $w_i=0.25$ for $i\in\{1,4\}$. \Gls{IQTS}: $m=4$, $n=15$, $\kappa=50$. \Gls{QAOA}: $p=3$. \emph{SV1}: \num{1024} shots per circuit. Total number of instances: \num{20}.
\end{itemizedense}

In all experiments, we choose a global denominator $R := 10$ for the relative values, \cref{eqn:vapprox}, and set the corresponding part-dependent numerators to their optimal values
\begin{align}
     \hat{P}_i := \argmin_{P_i \in \mathbb{Z}} \left| v_i - \frac{P_i}{Q} \right|
\end{align}
for all $i \in I$. This choice is particularly independent of the scalarization weights $w$ and leads to a maximum error of $\max_{i \in I} \epsilon_i = \num{0.051}$.

All algorithms are implemented in Python using common libraries. Pennylane~\cite{pennylane} is used to realize quantum circuits.

Quantum circuits for \gls{QAOA} are either executed on the IonQ device \emph{Aria-1} or simulated with the AWS Braket universal state vector simulator \emph{SV1}. Classical algorithms also run on AWS. For \gls{IQTS}, we use 4 instances of AWS EC2 in total, each with 4 (total 16) vCPUs (3.6 GHz \emph{Intel Xeon}) and 8 GiB of memory. For \gls{HBS}, we use 22 instances of AWS EC2 in total with each 36 (total 792) vCPUs (2.9 GHz \emph{Intel Xeon E5-2666}) and 60 GiB of memory.

\section{Additional Results} \label{supp:exp}

\begin{figure}[!htb]
	\centering
    \includegraphics[scale=.98]{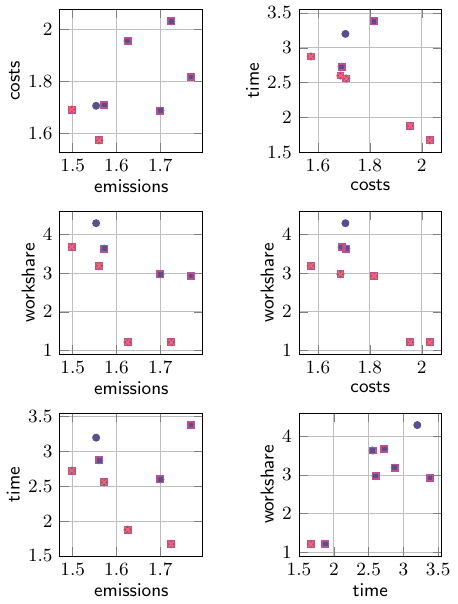}
	\caption{Solutions from \cref{exp:6} with the same symbols as in \cref{fig:pareto-e1}.}
	\label{fig:pareto-e6}
\end{figure}

\begin{figure}[!htb]
	\centering
    \includegraphics[scale=.98]{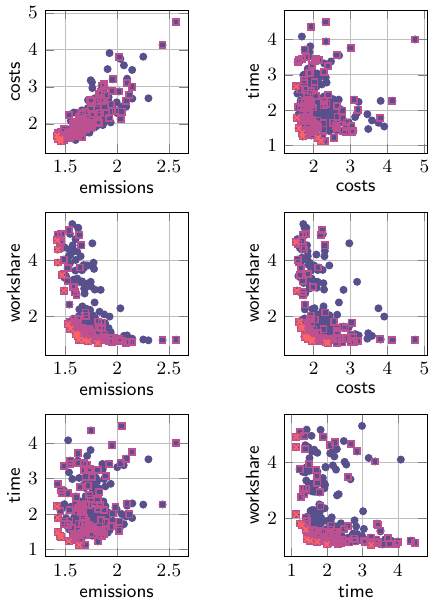}
	\caption{Solutions from \cref{exp:2} with the same symbols as in \cref{fig:pareto-e1}.}
	\label{fig:pareto-e2}
\end{figure}

\begin{figure}[!htb]
	\centering
    \includegraphics[scale=.98]{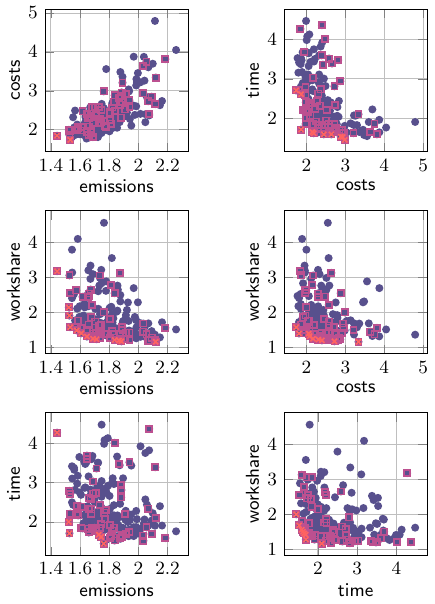}
	\caption{Solutions from \cref{exp:4} with the same symbols as in \cref{fig:pareto-e1}.}
	\label{fig:pareto-e4}
\end{figure}

\begin{figure}[!htb]
	\centering
    \includegraphics[scale=.98]{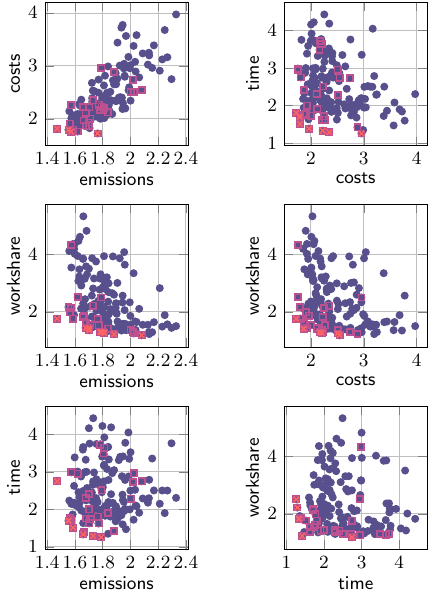}
	\caption{Solutions from \cref{exp:5} with the same symbols as in \cref{fig:pareto-e1}.}
	\label{fig:pareto-e5}
\end{figure}

\begin{figure}[!htb]
	\centering
    \includegraphics[scale=.98]{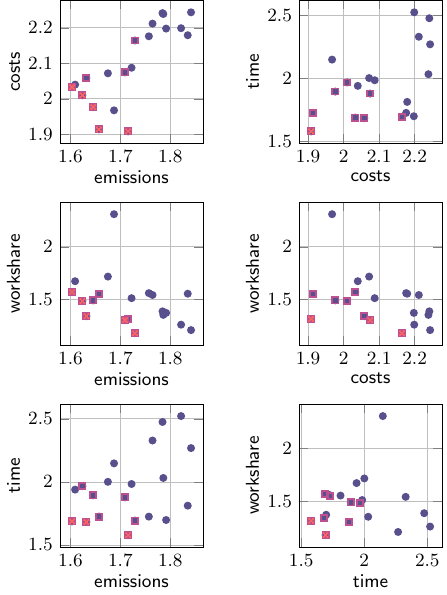}
	\caption{Solutions from \cref{exp:3} with the same symbols as in \cref{fig:pareto-e1}.}
	\label{fig:pareto-e3}
\end{figure}

\begin{figure}[!htb]
	\centering
    \includegraphics[scale=.98]{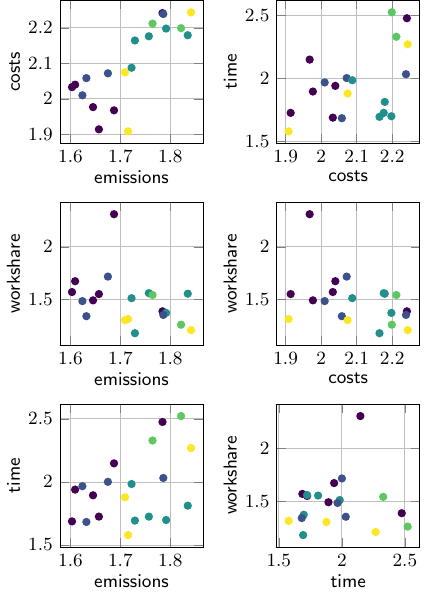}
	\caption{Solutions from \cref{exp:3} with the same symbols as in \cref{fig:alphas}.}
	\label{fig:alphas-supp}
\end{figure}

In the following, we provide additional results for \cref{sec:experiments}. Firstly, we present in \cref{fig:pareto-e6,fig:pareto-e2,fig:pareto-e4,fig:pareto-e5,fig:pareto-e3} the Pareto frontiers for \cref{exp:6,exp:2,exp:4,exp:5,exp:3} in analogy to \cref{fig:pareto-e1} for \cref{exp:1}. In each projection, we highlight the Pareto-optimal points in the four-dimensional \gls{KPI} space as well as the optimal points for this specific projection (when ignoring the other dimensions). Note that all presented solutions are feasible. Secondly, we present the remaining projections from \cref{exp:3} in \cref{fig:alphas-supp} as a supplement to \cref{fig:alphas}.

\end{appendices}

\end{document}